\documentclass[%
onecolumn, a4paper
]{mpi2015-cscpreprint}

\usepackage[american]{babel}

\usepackage{graphicx}

\usepackage{amssymb}
\usepackage{amsthm}

\usepackage{amsmath}
\usepackage{subfig}

\usepackage{cite}
\usepackage{threeparttable}
\usepackage{multirow}
\usepackage{float}
\usepackage{enumerate}
\usepackage[none]{hyphenat}
\usepackage[normalem]{ulem}
\usepackage{color}
\usepackage{xcolor}
\usepackage{mathtools}

\usepackage{tikz}
\usepackage{pgfplots}
\pgfplotsset{compat=newest}
\usetikzlibrary{positioning}

\usepackage{algorithm}
\usepackage{algorithmic}

\newtheorem{remark}{Remark}

\newlength\fheight
\newlength\fwidth
\begin{document}
	
%
\title{Fast \emph{A Posteriori} State Error Estimation for Reliable Frequency Sweeping in Microwave Circuits via the Reduced-Basis Method}
\author[$\dagger$]{Valent\'in de la Rubia}
\affil[$\dagger$]{Departamento de Matem\'{a}tica Aplicada a las {TIC}, ETSI de Telecomunicaci\'{o}n, Universidad Polit\'{e}cnica de Madrid, 28040 Madrid, Spain.\authorcr
	\email{valentin.delarubia@upm.es}, \orcid{0000-0002-2894-6813}}
\author[$\ast$]{Sridhar Chellappa}
\affil[$\ast$]{Department of Computational Methods in Systems and Control Theory, Max Planck Institute for Dynamics of Complex Technical Systems, 39106 Magdeburg, Germany.\authorcr
	\email{chellappa@mpi-magdeburg.mpg.de}, \orcid{0000-0002-7288-3880} \authorcr
	\email{feng@mpi-magdeburg.mpg.de}, \orcid{0000-0002-1885-3269} \authorcr
	\email{benner@mpi-magdeburg.mpg.de}, \orcid{0000-0003-3362-4103}}

\author[$\ast$]{Lihong Feng}

\author[$\ast$]{Peter Benner}

\shorttitle{Fast State Error Estimation for MOR of Microwave Circuits}
\shortauthor{V. de la Rubia, S. Chellappa, L. Feng, P. Benner}
\shortdate{}

\keywords{Computer aided engineering, design automation, error analysis, finite element method, Galerkin method, microwave circuits, reduced basis method, reduced order modeling.}

\msc{37M05, 65M60, 65R20, 68U07, 78A50}

\abstract{%
	We develop a compact, reliable model order reduction approach for fast frequency sweeps in microwave circuits by means of the reduced-basis method. Contrary to what has been previously done, special emphasis is placed on certifying the accuracy of the reduced-order model with respect to the original full-order model in an effective and efficient way. Previous works on model order reduction accuracy certification rely on costly \emph{a posteriori} error estimators, which typically require expensive \emph{inf-sup} constant evaluations of the underlying full-order model. This scenario is often too time-consuming and unaffordable in electromagnetic applications. As a result, less expensive and heuristic error estimators are commonly used instead. Very often, one is interested in knowing about the full state vector, instead of just some output quantities derived from the full state. Therefore, error estimators for the full state vector become relevant. In this work, we detail the frequency behavior of both the electric field and the state error when an approximation to the electric field solution is carried out. Both field quantities share the same frequency behavior. Based on this observation, we focus on the efficient estimation of the electric field state error and propose a fast evaluation of the reduced-order model state error in the frequency band of analysis, minimizing the number of full-order model evaluations. This methodology is of paramount importance to carry out a reliable fast frequency sweep in microwave circuits. Finally, real-life applications will illustrate the capabilities and efficiency of the proposed approach.
}
 \novelty{We propose a fast method for \emph{a posteriori} state error estimation for reduced-order modeling of microwave circuits. Our approach leverages a similarity in the frequency domain behavior of the electric field and the state error, which enables us to minimize the number of full-order model evaluations during error estimation. The proposed fast state error estimator is used in a greedy algorithm to adaptively build the reduced-order model. The efficiency and reliability of our approach are demonstrated via its application to four challenging models of microwave circuits.}

\maketitle
\section{Introduction}
\label{Sec-Introduction}
Microwave engineering relies on time-consuming electromagnetic simulations to carry out robust electrical designs. The electromagnetic complexity in microwave devices is such that \emph{only} detailed, full-wave simulations, i.e., solving Maxwell's equations directly, can guide an engineer in pursuing the target electrical design. As a result, most of a microwave engineer's working time is spent on waiting for an electromagnetic simulation to conclude, which will assist him in taking an action to meet the specifications in an electrical design. Current industrial needs have long since recognized that one can no longer afford this design methodology. Different efforts in computational electromagnetics (CEM) community have been carried out to speed up this costly process, and most of them follow the model order reduction (MOR) philosophy  \cite{Edlinger2014APosteriori,nicolini2019model,hochman2014reduced,Rewienski2016greedy,codecasa2019exploiting,ChewTAP2014GeneralizedModal,xue2020rapid,mrozowski2020,szypulski2020SSMMM}.

A reduced-order model (ROM) implies replacing a rather complex physical model by a much simpler mathematical one that still maintains certain physical aspects of the original model over a parameter domain. The computational complexity of the ROM should be insignificant in comparison to the high computational cost of the original full-order model (FOM). MOR has demonstrated its robustness in reducing the complexity of parametric systems \cite{morQuaMN16,morHesRS16,morBenGW15}. However, the accuracy of the ROM is sometimes not guaranteed due to lack of low-cost and computable error estimators. Although the ROM may be valid for a certain parameter range, its validity over the entire parameter domain is not guaranteed. As a result, it can not be used as a reliable surrogate of the original FOM. This lack of providing accuracy guarantees precludes the ROM from being used for industrial applications, where no \emph{a priori} knowledge of the parameter range may be available. This is the worst case scenario for MOR. It is quite usual that the proposition of corresponding error estimation lags behind new MOR algorithms. To remedy this, a great effort has been carried out in certifying the accuracy of the ROM, where computationally expensive error estimation may be allowed. This is the case for \emph{inf-sup} constant-based error estimators \cite{morHesB13,hess2015estimating,Edlinger2015CertifiedDualCorrected,morSchWH18}. The residual norm divided by this costly \emph{inf-sup} constant \cite{GarRM17} bounds the state error. As already stated in \cite{morFenB19}, keeping the \emph{inf-sup} constant in the denominator of the error estimator causes potential risk for many problems with small \emph{inf-sup} constants. This is quite common in microwave circuits, where resonances show up in the frequency band of analysis, dropping the \emph{inf-sup} constant down to zero \cite{delaRubia2018CRBM,GarRM17}. Different strategies for \emph{a posteriori} error estimation should be considered, reducing its computational cost to the same order of the ROM, if possible. Residual norm-based error estimation can be carried out without effort and this has been often used as a heuristic error estimator~\cite{delaRubia2018CRBM,de2009reliable, Vouvakis2011FastFrequency, delaRubia2014Reliable,  morRewLM15, Edlinger2015ANewMethod, Edlinger2017finite, fotyga2018reliable, MonjedelaRubia2020EFIE}. Going back to the state error estimation, recent works have focused on avoiding the \emph{inf-sup} constant evaluation \cite{morSemZP18, feng2020InfSupConstantFree}. There, additional dual or residual systems are solved to obtain the error estimators, overcoming any time-consuming \emph{inf-sup} constant calculation. However, despite the fact that they avoid computing the expensive \emph{inf-sup} constant, both approaches need to solve additional dual or residual systems, respectively.   

In this work, we aim to further reduce the computational costs of the error estimator proposed in \cite{feng2020InfSupConstantFree}, which was shown to be more efficient than that in \cite{morSemZP18}. We study the frequency behavior of both the electric field and the state error when an approximation to the electric field solution is carried out, and detail a Fourier series representation in both cases. Both field quantities share the same orthogonal Fourier series representation basis in frequency-parameter systems. In comparison to what has been previously done for general parametric systems, where a residual system needs to be solved independently, we focus on the efficient determination of the electric field state error and propose a fast evaluation of the ROM state error in the frequency band of analysis, minimizing the number of FOM evaluations which plays a central role in determining the efficiency in MOR. This methodology is of paramount importance to carry out a reliable fast frequency sweep in microwave circuits.

This paper is organized as follows. In Section \ref{Sec-ProblemStatement} we review the time-harmonic Maxwell's equations in variational form, solve for the electromagnetic field in order to show its frequency behavior and detail the \emph{inf-sup} constant-based standard error analysis. Section \ref{Sec-InfSupConstantFree} deals with the proposed state error estimation in frequency-parameter systems avoiding the \emph{inf-sup} constant. Numerical simulations in Section \ref{Sec-NumericalResults} show the performance of the proposed approach for reliable fast frequency sweeps in electromagnetics. Real-life microwave circuits illustrate the capabilities and accuracy of the proposed methodology. Finally, in Section \ref{Sec-Conclusions}, we provide conclusions.
\section{Problem Statement}
\label{Sec-ProblemStatement}
The electromagnetic phenomena in a given device are described by Maxwell's equations. Applying the Fourier transform to these, the fields in the transform domain $i\omega$ can be found. They are
\begin{subequations}
	\label{eq:Sec-ProblemStatement-MaxwellSystem}
	\begin{align}
	\nabla \times \mathbf{E} &= -i \omega \mu \mathbf{H} \text{ in } \Omega,\\
	\nabla \times \mathbf{H} &= i \omega \varepsilon \mathbf{E} \text{ in } \Omega,\\
	\label{eq:Sec-ProblemStatement-MaxwellSystem-PECBoundaryCondition}
	\mathbf{n} \times \mathbf{E} &= \mathbf{0} \text{ on } \Gamma_\text{PEC},\\
	\mathbf{n} \times \mathbf{H} &= \mathbf{0} \text{ on } \Gamma_\text{PMC},\\
	\mathbf{n} \times \mathbf{H} &= \mathbf{J} \text{ on } \Gamma \text{,}
	\end{align}
\end{subequations}
where $\Omega \subset \mathbb{R}^{3}$ is a source-free, sufficiently smooth bounded domain, $\mathbf{n}$ is the unit outward normal vector on the boundary $\partial \Omega$ of $\Omega$. The boundary is divided into perfect electric conductor (PEC), perfect magnetic conductor (PMC) and ports, i.e., $\partial \Omega =\Gamma_\text{PEC} \cup \Gamma_\text{PMC}\cup \Gamma $. $\mathbf{E}$ and $\mathbf{H}$ are the electric and magnetic  fields, $\varepsilon $ and $\mu $ are, respectively, the permittivity and permeability of the medium, which is assumed to be lossless, and the tangential field $\mathbf{J}$ is the excitation current at the ports. Time-harmonic Maxwell's equations can be written in a classical weak formulation over an appropriate admissible function space $\mathcal{H}$, viz.
\begin{equation}
\label{eq:Sec-ProblemStatement-WeakForm}
\begin{aligned}
\text{find}~&\mathbf{E}\in \mathcal{H}~\text{such that} \\ 
&a(\mathbf{E},\mathbf{v})=f(\mathbf{v})~\forall \mathbf{v}\in \mathcal{H}\text{.}
\end{aligned}
\end{equation}
The bilinear form is
\begin{equation}
\label{eq:Sec-ProblemStatement-BilinearForm}
a(\mathbf{E}, \mathbf{v}) = \int\limits_{\Omega}\left( \frac{1}{\mu } \nabla \times \mathbf{E} \cdot \nabla \times \mathbf{v} - \omega ^{2} \varepsilon \mathbf{E} \cdot \mathbf{v}\right) dx \text{,}
\end{equation}
and the linear form
\begin{equation}
\label{eq:Sec-ProblemStatement-LinearForm}
f(\mathbf{v})=i\omega \int\limits_{\partial \Omega} \mathbf{J} \cdot \mathbf{v}~ds = i \omega \int\limits_{\Gamma} \mathbf{J} \cdot \mathbf{v}~ds \text{.}
\end{equation}
Here, the admissible space $\mathcal{H}$ is a subspace of the Hilbert space $H(curl,\Omega )$ defined by:
\begin{equation}
\label{eq:Sec-ProblemStatement-Hcurl}
H(curl,\Omega )=\left\{ \mathbf{u}\in L^{2}(\Omega, \mathbb{C}^{3})~|~\nabla \times \mathbf{u}\in L^{2}(\Omega, \mathbb{C}^{3})\right\} \text{,}
\end{equation}
since $\mathcal{H}$ should take the boundary condition \eqref{eq:Sec-ProblemStatement-MaxwellSystem-PECBoundaryCondition} into account, namely,
\begin{equation}
\label{eq:Sec-ProblemStatement-H}
\mathcal{H}=\left\{ \mathbf{u}\in H(curl, \Omega )~|~\mathbf{n} \times \mathbf{u} = \mathbf{0}~\text{on }\Gamma_\text{PEC}\right\} \text{.}
\end{equation}

Let us refer to the trace spaces, namely,
\begin{equation}
\label{eq:Sec-ProblemStatement-TraceSpaces}
\begin{aligned}
H^{-\frac{1}{2}}(div, \partial \Omega)&= \{ \mathbf{n} \times \mathbf{u}~\text{on}~\partial \Omega ~|~\mathbf{u} \in H(curl,\Omega) \} \\
H^{-\frac{1}{2}}(curl, \partial \Omega)&= \{ \mathbf{n} \times \mathbf{u} \times \mathbf{n}~\text{on}~\partial \Omega ~|~\mathbf{u} \in H(curl,\Omega) \} \text{,}
\end{aligned}
\end{equation}
and point out that they are dual to each other with the following duality pairing
\begin{equation}
\label{eq:Sec-ProblemStatement-DualityPairing}
\langle \mathbf{u},\mathbf{v} \rangle = \int\limits_{\partial \Omega} \mathbf{u} \cdot \mathbf{v} ~ds
\end{equation}
$\mathbf{u} \in H^{-\frac{1}{2}}(div,\partial \Omega)$ and $\mathbf{v} \in H^{-\frac{1}{2}}(curl,\partial \Omega)$. It is now apparent that the excitation current $\mathbf{J}$ belongs to $H^{-\frac{1}{2}}(div, \partial \Omega)$. We refer to \cite{GiraultRaviart,Monk} for a through explanation for all these spaces.
\subsection{Field Frequency Dependency in Electromagnetics}
Following \cite{delaRubia2018CRBM,Kirsch}, where some frequency structure is shown in the solution to the variational problem \eqref{eq:Sec-ProblemStatement-WeakForm}, we introduce the Helmholtz decomposition
\begin{equation}
\label{eq:Sec-ProblemStatement-HelmholtzDecomposition}
\mathcal{H}=\mathcal{H}(curl0, \Omega) \oplus \mathcal{V} \text{,}
\end{equation}
where
\begin{subequations}
	\label{eq:Sec-ProblemStatement-HelmholtzSpaces}
	\begin{align}
	\label{eq:Sec-ProblemStatement-HelmholtzSpaces-Hcurl0}
	&\mathcal{H}(curl0,\Omega) = \{ \mathbf{u} \in \mathcal{H}~|~ \nabla \times \mathbf{u} = \mathbf{0} \},\\
	\label{eq:Sec-ProblemStatement-HelmholtzSpaces-V}
	&\mathcal{V} = \{ \mathbf{u} \in \mathcal{H}~|~(\varepsilon \mathbf{u}, \mathbf{v})_{L^{2}(\Omega)} = 0~\forall \mathbf{v} \in \mathcal{H}(curl0,\Omega) \} \text{.}
	\end{align}
\end{subequations}
$(\cdot, \cdot)_{L^2(\Omega)}$ is the inner product in $L^{2}(\Omega, \mathbb{C}^{3})$. $\mathcal{H}(curl0, \Omega)$ denotes the nullspace of the curl operator while $\mathcal{V}$ stands for its orthogonal complement within the solution space $\mathcal{H}$ in the following inner product
\begin{equation}
\label{eq:Sec-ProblemStatement-muepsInnerProduct}
(\mathbf{u}, \mathbf{v})_{\mu, \varepsilon} = ( \frac{1}{\mu} \nabla \times \mathbf{u}, \nabla \times \mathbf{v})_{L^{2}(\Omega)}+(\varepsilon \mathbf{u}, \mathbf{v})_{L^{2}(\Omega)} \text{.}
\end{equation}
It should be noted that both $\mathcal{H}(curl0, \Omega)$ and $\mathcal{V}$ spaces satisfy the PEC boundary condition on $\Gamma_\text{PEC}$.

The variational problem \eqref{eq:Sec-ProblemStatement-WeakForm} can be solved by using the splitting $\mathbf{E} = \mathbf{E}_0 + \mathbf{e}$, $\mathbf{E}_0 \in~\mathcal{H}(curl0, \Omega)$, $\mathbf{e} \in \mathcal{V}$. We refer to \cite{delaRubia2018CRBM,Kirsch} for details. As a result, we can make the dependence of the solution to time-harmonic Maxwell's equations on frequency explicit, cf. \cite{Kurokawa,Conciauro}, i.e.,
\begin{equation}
\label{eq:Sec-ProblemStatement-SolutionMaxwellFrequencyDependency}
\begin{aligned}
&\text{if } \omega^2 \neq \omega_n^2\text{,} ~ \mathbf{E} = \mathbf{E}_0 + \mathbf{e} = \frac{1}{i\omega} \mathbf{F}_0 + i \omega \sum \limits_{n=1}^{\infty} \frac{A_n}{1 - \frac{\omega^2}{\omega_n^2} } \mathbf{e}_n,\\
&\text{if } \omega^2 = \omega_n^2\text{,}~\mathbf{E} = \mathbf{E}_0 + \mathbf{e} = \\
&= \frac{1}{i\omega} \mathbf{F}_0 + i \omega \sum \limits_{ \omega_n^2 = \omega^2 } a_n \mathbf{e}_n + i \omega \sum \limits_{\omega_n^2 \neq \omega^2} \frac{A_n}{ 1 - \frac{\omega^2}{\omega_n^2} } \mathbf{e}_n \text{.}
\end{aligned}
\end{equation}
$\mathbf{F}_0 \in \mathcal{H}(curl0, \Omega)$ is related to the Riesz representative for the electric field in statics. The set of eigenmodes $\{\mathbf{e}_n~|~n \in~\mathbb{N}\} \subset \mathcal{V}$ stands for the resonant modes in electrodynamics, along with their corresponding resonant frequencies $\omega_n \in \mathbb{R}$, and forms a complete orthonormal system in $\mathcal{V}$ with respect to the inner product \eqref{eq:Sec-ProblemStatement-muepsInnerProduct} \cite{Kirsch}. It should be pointed out that $\mathcal{H}(curl0, \Omega)$ is orthogonal to $\mathcal{V}$ with respect to the same inner product \eqref{eq:Sec-ProblemStatement-muepsInnerProduct}. Getting to our point, \eqref{eq:Sec-ProblemStatement-SolutionMaxwellFrequencyDependency} details an orthogonal representation, i.e., a Fourier series for the electric field where the frequency dependence is explicit. Further, $A_n$ are coupling coefficients for the excitation current $\mathbf{J}$ to its corresponding resonant mode $\mathbf{e}_n$ and are determined by 
\begin{equation}
\label{eq:Sec-ProblemStatement-ExcitationCouplingCoefficient}
A_n = \langle \mathbf{J}, \mathbf{n} \times \overline{\mathbf{e}}_n \times \mathbf{n} \rangle \text{.}
\end{equation}
$\overline{\mathbf{e}}_n$ stands for the complex conjugate of $\mathbf{e}_n$. Finally, $a_n$ are arbitrary coefficients since the electric field is not unique at resonance.
\subsection{Parametric Variational Problem and Standard \emph{A Posteriori} Error Analysis}
\label{Sec-ProblemStatement-Subsec-StandarAPosterioriError}
Taking frequency ($\omega$) as a parameter, the weak formulation for time-harmonic Maxwell's equations \eqref{eq:Sec-ProblemStatement-WeakForm} turns into the following parametric variational problem:
\begin{equation}
\label{eq:Sec-ProblemStatement-Subsec-StandarAPosterioriError-ParametricWeakForm}
\begin{aligned}
\text{find}~&\mathbf{E}(\omega) \in \mathcal{H}~\text{such that} \\ 
&a(\mathbf{E}(\omega), \mathbf{v}; \omega) = f(\mathbf{v}; \omega)~\forall \mathbf{v} \in \mathcal{H},~\forall \omega \in \mathcal{B} \text{.}
\end{aligned}
\end{equation}
where $\mathcal{B}:=[\omega_\text{min},\omega_\text{max}]\subset \mathbb{R}$ is the frequency band of interest, and the frequency-parameter bilinear and linear forms $a(\cdot,\cdot;\omega)$ and $f(\cdot;\omega)$ are already defined in \eqref{eq:Sec-ProblemStatement-BilinearForm} and \eqref{eq:Sec-ProblemStatement-LinearForm}, respectively. The well-posedness of the parametric problem \eqref{eq:Sec-ProblemStatement-Subsec-StandarAPosterioriError-ParametricWeakForm} relies on the behavior of the so-called \emph{inf-sup} constant $\beta(\omega)$ as a function of frequency:
\begin{equation}
\label{eq:Sec-ProblemStatement-Subsec-StandarAPosterioriError-InfSupConstant}
\beta (\omega ) = \adjustlimits\inf_{\mathbf{u} \in \mathcal{H}} \sup_{\mathbf{v} \in \mathcal{H}}\, \frac{ |a(\mathbf{u}, \mathbf{v}; \omega)| }{ \| \mathbf{u} \|_\mathcal{H} \| \mathbf{v} \|_\mathcal{H} }.
\end{equation}
For all $\omega \in \mathcal{B}$, $\beta (\omega )\geq \beta _{0} > 0$ ensures the well-posedness and uniqueness in the variational problem \eqref{eq:Sec-ProblemStatement-Subsec-StandarAPosterioriError-ParametricWeakForm} \cite{hess2015estimating}. 

This result gives rise to a standard \emph{a posteriori} error analysis.  Provided an approximate solution $\mathbf{\tilde E}(\omega) \in \mathcal{H}$ to the variational problem \eqref{eq:Sec-ProblemStatement-Subsec-StandarAPosterioriError-ParametricWeakForm} is found, the error in the approximation $\| \mathbf{E}(\omega) - ~\mathbf{\tilde E}(\omega)\|_\mathcal{H}$ can be bounded using the \emph{inf-sup} constant. Indeed, \eqref{eq:Sec-ProblemStatement-Subsec-StandarAPosterioriError-InfSupConstant} can be rewritten as follows:
\begin{equation}
\label{eq:Sec-ProblemStatement-Subsec-StandarAPosterioriError-InfSupConstant2}
\sup_{\mathbf{v} \in \mathcal{H}} \frac{|a(\mathbf{u}, \mathbf{v};\omega)|}{ \|\mathbf{u}\|_\mathcal{H} \|\mathbf{v}\|_\mathcal{H}}>\beta(\omega), \forall \mathbf{u} \in \mathcal{H}.  
\end{equation}
In particular, this inequality still holds when replacing $\mathbf{u} \in \mathcal{H}$ by the field $\mathbf{E}(\omega)-\mathbf{\tilde E}(\omega) \in \mathcal{H}$, which gives rise to an upper bound for the approximation error, namely, 
\begin{equation}
\label{eq:Sec-ProblemStatement-Subsec-StandarAPosterioriError-UpperBoundError}
\| \mathbf{E}(\omega )- \mathbf{ \tilde E}(\omega) \|_\mathcal{H} < \frac{ 1 }{ \beta(\omega)} \sup_{ \mathbf{v} \in \mathcal{H} } \frac{ |a(\mathbf{E}(\omega )-\mathbf{\tilde E}(\omega), \mathbf{v}; \omega)| }{ \| \mathbf{v} \|_\mathcal{H} }.
\end{equation}
However, this error bound not only involves the computation of the norm of the residual functional 
\begin{equation}
\label{eq:Sec-ProblemStatement-Subsec-StandarAPosterioriError-Residual}
\begin{aligned}
r(\mathbf{\tilde E}(\omega), \mathbf{v}; \omega) :=& f(\mathbf{v}; \omega )-a(\mathbf{\tilde E}(\omega), \mathbf{v}; \omega) \\
=& a(\mathbf{E}(\omega)-\mathbf{\tilde E}(\omega), \mathbf{v}; \omega),\forall \mathbf{v} \in \mathcal{H} \text{,}
\end{aligned}
\end{equation}
which can be determined in an efficient way as a function of frequency \cite{de2009reliable,Vouvakis2011FastFrequency,morHesB13,Edlinger2015ANewMethod,morRewLM15}, but also the determination of the \emph{inf-sup} constant  throughout the frequency band of interest $\mathcal{B}$, which can be time-consuming \cite{hess2015estimating,Edlinger2015CertifiedDualCorrected,GarRM17}.

Furthermore, in microwave engineering, resonances appearing in $\mathcal{B}$ are responsible for the target electrical response. As a result, resonant modes arise and the uniqueness of the solution is no longer valid in the band of interest $\mathcal{B}$. The \emph{inf-sup} constant vanishes at the resonance frequencies, giving rise to a near-infinity upper bound for the error in \eqref{eq:Sec-ProblemStatement-Subsec-StandarAPosterioriError-UpperBoundError} nearby resonances. In addition, the above error estimation leads to unacceptable overestimation of the error even for well-conditioned problems~\cite{morSchWH18}. Given no better choices, the norm of the residual \eqref{eq:Sec-ProblemStatement-Subsec-StandarAPosterioriError-Residual}, which can be straightforwardly computed, has been used as a \emph{heuristic} error estimator \cite{Rewienski2016greedy,delaRubia2014Reliable,morRewLM15,Edlinger2017finite,fotyga2018reliable,MonjedelaRubia2020EFIE,Kouki}
\section{State Error Estimation Avoiding the \emph{Inf-Sup}-Constant}
\label{Sec-InfSupConstantFree}
The previous section has shown the main role the \emph{inf-sup} constant plays in \emph{a posteriori} error estimation, as well as the incapability of \emph{inf-sup} constant-based error estimators to provide a tight error bound nearby resonant frequencies. Unfortunately, the norm of the residual cannot provide a sharp error estimation at or nearby resonance frequencies either. We shall elaborate on this later in Section \ref{Sec-Conclusions}. As a result, we are in need of more efficient state error estimators to certify the accuracy of the approximate field solution to the frequency-parameter variational problem \eqref{eq:Sec-ProblemStatement-Subsec-StandarAPosterioriError-ParametricWeakForm}, even in the presence of resonances. Recent efforts have moved towards this goal \cite{morSemZP18,feng2020InfSupConstantFree}. There, instead of computing the \emph{inf-sup} constant, additional dual or residual systems need to be solved to obtain the state error estimator. These additional systems constitute an extra computational effort to certify the accuracy of the approximate field solution. In this work, we focus on fast \emph{a posteriori} state error estimator computation taking advantage of the frequency dependency in the field solution \eqref{eq:Sec-ProblemStatement-SolutionMaxwellFrequencyDependency} for the frequency-parameter problem \eqref{eq:Sec-ProblemStatement-Subsec-StandarAPosterioriError-ParametricWeakForm}.
\subsection{Field Error Frequency Dependency in Electromagnetics}
Given an approximate solution $\mathbf{\tilde E}(\omega) \in \mathcal{H}$ to \eqref{eq:Sec-ProblemStatement-Subsec-StandarAPosterioriError-ParametricWeakForm}, we can study the Fourier series representation of the state error 
\begin{equation}
\label{eq:Sec-InfSupConstantFree-StateError}
	\boldsymbol{\epsilon}(\omega) := \mathbf{E}(\omega) - \mathbf{\tilde E}(\omega) \in \mathcal{H} \text{,}
\end{equation}
making its frequency dependency explicit. The state error \eqref{eq:Sec-InfSupConstantFree-StateError} satisfies the frequency-parameter variational problem,
\begin{equation}
\label{eq:Sec-InfSupConstantFree-ParametricWeakFormStateError}
\begin{aligned}
\text{find}~&\boldsymbol{\epsilon}(\omega) \in \mathcal{H}~\text{such that} \\ 
&a(\boldsymbol{\epsilon}(\omega), \mathbf{v}; \omega) = f^{\boldsymbol{\epsilon}}(\mathbf{v}; \omega)~\forall \mathbf{v} \in \mathcal{H},~\forall \omega \in \mathcal{B} \text{.}
\end{aligned}
\end{equation}
The frequency-parameter bilinear form $a(\cdot, \cdot; \omega)$ is already defined in \eqref{eq:Sec-ProblemStatement-BilinearForm}. The frequency-parameter linear form $f^{\boldsymbol{\epsilon}}(\cdot; \omega)$ is the residual functional $r(\mathbf{\tilde E}(\omega), \cdot; \omega) $ detailed in \eqref{eq:Sec-ProblemStatement-Subsec-StandarAPosterioriError-Residual},
\begin{equation}
\label{eq:Sec-InfSupConstantFree-LinearForm}
f^{\boldsymbol{\epsilon}}(\mathbf{v}) = f(\mathbf{v}; \omega )-a(\mathbf{\tilde E}(\omega), \mathbf{v}; \omega) = i \omega \int\limits_{\partial \Omega} \mathbf{J}^{\boldsymbol{\epsilon}} \cdot \mathbf{v}~ds \text{,}
\end{equation}
which can be identified as a residual error current $\mathbf{J}^{\boldsymbol{\epsilon}} \in ~ H^{-\frac{1}{2}}(div, \partial \Omega)$. By an analogous reasoning as the one in Section \ref{Sec-ProblemStatement}, we get
\begin{equation}
\label{eq:Sec-InfSupConstantFree-SolutionErrorMaxwellFrequencyDependency}
\begin{aligned}
\text{if } \omega^2 &\neq \omega_n^2\text{,} ~ \boldsymbol{\epsilon}(\omega) = \frac{1}{i\omega} \mathbf{F}^{\boldsymbol{\epsilon}}_0 + i \omega \sum \limits_{n=1}^{\infty} \frac{A^{\boldsymbol{\epsilon}}_n}{1 - \frac{\omega^2}{\omega_n^2} } \mathbf{e}_n,\\
\text{if } \omega^2 &= \omega_n^2\text{,} \\
\boldsymbol{\epsilon}(\omega) &= \frac{1}{i\omega} \mathbf{F}^{\boldsymbol{\epsilon}}_0 + i \omega \sum \limits_{ \omega_n^2 = \omega^2 } a^{\boldsymbol{\epsilon}}_n \mathbf{e}_n + i \omega \sum \limits_{\omega_n^2 \neq \omega^2} \frac{A^{\boldsymbol{\epsilon}}_n}{ 1 - \frac{\omega^2}{\omega_n^2} } \mathbf{e}_n \text{.}
\end{aligned}
\end{equation}
$\mathbf{F}^{\boldsymbol{\epsilon}}_0 \in \mathcal{H}(curl0, \Omega)$ is related to the Riesz representative for the stationary error field. The same set of eigenmodes $\{\mathbf{e}_n~|~n \in~\mathbb{N}\} \subset \mathcal{V}$ along with their corresponding resonant frequencies $\omega_n$ as in \eqref{eq:Sec-ProblemStatement-SolutionMaxwellFrequencyDependency} can be used. $A^{\boldsymbol{\epsilon}}_n$ are coupling coefficients for the residual error current $\mathbf{J}^{\boldsymbol{\epsilon}}$ to the corresponding resonant mode $\mathbf{e}_n$, namely,
\begin{equation}
\label{eq:Sec-InfSupConstantFree-ErrorExcitationCouplingCoefficient}
A^{\boldsymbol{\epsilon}}_n = \langle \mathbf{J}^{\boldsymbol{\epsilon}}, \mathbf{n} \times \overline{\mathbf{e}}_n \times \mathbf{n} \rangle \text{.}
\end{equation}
In addition, $a^{\boldsymbol{\epsilon}}_n$ are arbitrary coefficients since there is no unique solution at resonance.

Having a closer look at equations \eqref{eq:Sec-ProblemStatement-SolutionMaxwellFrequencyDependency} and \eqref{eq:Sec-InfSupConstantFree-SolutionErrorMaxwellFrequencyDependency}, we can realize that the solutions to both original and residual variational problems share the same frequency-parameter behaviour and admit a similar Fourier series representation with the same frequency pattern. We may be then tempted to use the same representation basis to find an approximate solution to both original and residual variational problems \eqref{eq:Sec-ProblemStatement-Subsec-StandarAPosterioriError-ParametricWeakForm} and \eqref{eq:Sec-InfSupConstantFree-ParametricWeakFormStateError}. However, this should be done carefully to avoid stating that the state error $\boldsymbol{\epsilon}(\omega)$ is identically zero even in the situation where the approximate field $\mathbf{\tilde E}(\omega)$ may still be far away from the true solution $\mathbf{E}(\omega)$. We will get back to this point later in Section \ref{Sec-InfSupConstantFree-SubSec-InBandEigenmodes}. In order to avoid this embarrassing scenario, \cite{feng2020InfSupConstantFree} proposes to approximate these variational problems by applying different Galerkin projection spaces $\mathcal{H}_m, \mathcal{H}^{\boldsymbol{\epsilon}}_m:=\mathcal{H}_m+\mathcal{H}^\mathbf{r}_m \subset \mathcal{H}$, to the original and residual problems, respectively, giving rise to corresponding reduced systems which can be solved with ease, namely,
\begin{subequations}
\label{eq:Sec-InfSupConstantFree-GalerkinProjection}
\begin{align}
\label{eq:Sec-InfSupConstantFree-GalerkinProjectionField}
\text{find}~&\mathbf{\tilde E}(\omega) \in \mathcal{H}_m~\text{such that} \\ 
&a(\mathbf{\tilde E}(\omega), \mathbf{v}; \omega) = f(\mathbf{v}; \omega)~\forall \mathbf{v} \in \mathcal{H}_m,~\forall \omega \in \mathcal{B} \text{,} \nonumber \\
\label{eq:Sec-InfSupConstantFree-GalerkinProjectionError}
\text{and find}~&\boldsymbol{\tilde \epsilon}(\omega) \in \mathcal{H}^{\boldsymbol{\epsilon}}_m~\text{such that} \\ 
&a(\boldsymbol{\tilde \epsilon}(\omega), \mathbf{v}; \omega) = f^{\boldsymbol{\epsilon}}(\mathbf{v}; \omega)~\forall \mathbf{v} \in \mathcal{H}^{\boldsymbol{\epsilon}}_m,~\forall \omega \in \mathcal{B} \text{.} \nonumber
\end{align}
\end{subequations}
Algorithm \ref{alg:Sec-InfSupConstantFree-ROMErrorEstimator} adaptively builds the reduced-basis spaces up, i.e., the Galerkin projection spaces, in a greedy framework. As the number of iterations in this procedure increases, $\boldsymbol{\tilde \epsilon}(\omega)$ approximates the true state error $\boldsymbol{\epsilon}(\omega)$ better and better, $\boldsymbol{\tilde \epsilon}(\omega) \approx \boldsymbol{\epsilon}(\omega)$, and $\|\boldsymbol{\tilde \epsilon}(\omega)\|$ does perform as a sharp \emph{a posteriori} state error estimator. We refer to \cite{feng2020InfSupConstantFree} for the details. However, it should be pointed out that, in Algorithm \ref{alg:Sec-InfSupConstantFree-ROMErrorEstimator}, the dimension of the Galerkin projection space for the residual problem $\mathcal{H}^{\boldsymbol{\epsilon}}_m$ may double the dimension in the Galerkin projection space for the original problem $\mathcal{H}_m$ at each iteration. Further, distinct sets of parameters are arranged to solve for both the original and residual problems \eqref{eq:Sec-ProblemStatement-Subsec-StandarAPosterioriError-ParametricWeakForm} and \eqref{eq:Sec-InfSupConstantFree-ParametricWeakFormStateError}, which need to be done independently, in spite of the same dynamics being observed in both variational problems (see \eqref{eq:Sec-ProblemStatement-SolutionMaxwellFrequencyDependency} and \eqref{eq:Sec-InfSupConstantFree-SolutionErrorMaxwellFrequencyDependency}). This may turn this procedure rather time consuming. In this work, we focus on keeping a low computational effort, taking advantage of the observations in Sections \ref{Sec-ProblemStatement} and \ref{Sec-InfSupConstantFree} to carry out a reliable fast frequency sweep analysis.
\begin{algorithm}[t!]
	\caption{Adaptive construction of the Galerkin projection spaces $\mathcal{H}_m, \mathcal{H}^{\boldsymbol{\epsilon}}_m$ for state error estimation \cite{feng2020InfSupConstantFree}.}
	\label{alg:Sec-InfSupConstantFree-ROMErrorEstimator}
	\begin{algorithmic}[1]
		\REQUIRE Frequency band of interest $\mathcal{B}:=[\omega_\text{min}, \omega_\text{max} ]$, tolerance $\texttt{tol} > 0$ as the acceptable state error.
		\ENSURE $\mathcal{H}_m$ to ensure $\texttt{tol}$ state error in \eqref{eq:Sec-InfSupConstantFree-GalerkinProjection}.
		\STATE Initialize $\mathcal{H}_m = \{\mathbf{0}\}$, $\mathcal{H}^{\mathbf{r}}_m = \{\mathbf{0}\}$, $\xi = \texttt{tol} + 1$. Choose different samples $\omega^*$ and $\omega^*_{\boldsymbol{\epsilon}}$ randomly taken from $\mathcal{B}$.
		\WHILE{$\xi > \texttt{tol}$}
		\STATE Solve for $\mathbf{E}(\omega^*)$ in \eqref{eq:Sec-ProblemStatement-Subsec-StandarAPosterioriError-ParametricWeakForm}, orthonormalize and enrich $\mathcal{H}_m$: $\mathcal{H}_m = \mathcal{H}_m + \text{span}\{\mathbf{E}(\omega^*)\} $.
		\STATE Using $\mathcal{H}_m$, solve for the approximate field $\mathbf{\tilde E}(\omega)$ in \eqref{eq:Sec-InfSupConstantFree-GalerkinProjectionField}.
		\STATE Solve for $\boldsymbol{\epsilon}(\omega^*_{\boldsymbol{\epsilon}})$ in \eqref{eq:Sec-InfSupConstantFree-ParametricWeakFormStateError}, orthonormalize and enrich $\mathcal{H}^{\mathbf{r}}_m$: $\mathcal{H}^{\mathbf{r}}_m = \mathcal{H}^{\mathbf{r}}_m + \text{span}\{\boldsymbol{\epsilon}(\omega^*_{\boldsymbol{\epsilon}})\} $.
		\STATE Form $\mathcal{H}^{\boldsymbol{\epsilon}}_m$ and orthonormalize: $\mathcal{H}^{\boldsymbol{\epsilon}}_m = \mathcal{H}_m + \mathcal{H}^{\mathbf{r}}_m$.
		\STATE Using $\mathcal{H}^{\boldsymbol{\epsilon}}_m$, solve for the state error estimation $\boldsymbol{\tilde \epsilon}(\omega)$ in \eqref{eq:Sec-InfSupConstantFree-GalerkinProjectionError}.
		\STATE Choose the next sample $\omega^*$ from $\mathcal{B}$ as $$\omega^*=\arg\max\limits_{\omega\in \mathcal{B}}\|\boldsymbol{\tilde \epsilon}(\omega)\|_{\mathcal{H}}.$$
		\STATE Choose the next sample $\omega^*_{\boldsymbol{\epsilon}}$ from $\mathcal{B}$ following $$\omega^*_{\boldsymbol{\epsilon}}=\arg\max\limits_{\omega \in \mathcal{B}} \|r(\boldsymbol{\tilde \epsilon}(\omega), \cdot; \omega)\|_{\mathcal{H}^{\prime}}.$$ 
		Here $r(\boldsymbol{\tilde \epsilon}(\omega), \cdot; \omega)$ is the residual functional introduced by $\boldsymbol{\tilde \epsilon}(\omega)$ in \eqref{eq:Sec-InfSupConstantFree-ParametricWeakFormStateError}, namely, 
		$$r(\boldsymbol{\tilde \epsilon}(\omega), \mathbf{v}; \omega) := f^{\boldsymbol{\epsilon}}(\mathbf{v}; \omega )-a(\boldsymbol{\tilde \epsilon}(\omega), \mathbf{v}; \omega),\forall \mathbf{v} \in \mathcal{H}.$$
		\STATE $\xi = \|\boldsymbol{\tilde \epsilon}(\omega^*)\|_{\mathcal{H}}.$
		\label{step:Sec-InfSupConstantFree-Alg-ROMErrorEstimator-StoppingCriterion}	
		\ENDWHILE
		\STATE Use $\mathcal{H}_m$ to solve \eqref{eq:Sec-InfSupConstantFree-GalerkinProjectionField}.
	\end{algorithmic}
\end{algorithm}
%
\subsection{In-Band Eigenmodes in the Reduced-Basis Space}
\label{Sec-InfSupConstantFree-SubSec-InBandEigenmodes}
The works in \cite{szypulski2020SSMMM,delaRubia2018CRBM} suggest that a good approximation basis to the frequency-parameter problem \eqref{eq:Sec-ProblemStatement-Subsec-StandarAPosterioriError-ParametricWeakForm}, i.e., a good reduced-basis space, should include the resonant modes hit in the band of analysis $\mathcal{B}$. These resonant modes constitute the dominant basis representing the electric field $\mathbf{E}(\omega)$ in $\mathcal{B}$ (see \eqref{eq:Sec-ProblemStatement-SolutionMaxwellFrequencyDependency}). In this work, we use this basis not only for the original problem \eqref{eq:Sec-ProblemStatement-Subsec-StandarAPosterioriError-ParametricWeakForm}, but also for the residual problem \eqref{eq:Sec-InfSupConstantFree-ParametricWeakFormStateError}. For example,
\begin{subequations}
	\label{eq:Sec-InfSupConstantFree-SubSec-InBandEigenmodes-EigenBasis}
	\begin{align}
	\label{eq:Sec-InfSupConstantFree-SubSec-InBandEigenmodes-EigenBasisField}
	\mathcal{H}_m &= \text{span}\{\mathbf{e}_n~|~\omega_n \in~\mathcal{B}\}:=\mathcal{V}_\mathcal{B},\\
	\label{eq:Sec-InfSupConstantFree-SubSec-InBandEigenmodes-EigenBasisError}
	\mathcal{H}^{\boldsymbol{\epsilon}}_m &= \mathcal{H}_m = \mathcal{V}_\mathcal{B} \text{.}
	\end{align}
\end{subequations}
In this situation, we can find the solution for the reduced systems in \eqref{eq:Sec-InfSupConstantFree-GalerkinProjection} in closed form and get some further insights, viz.,
\begin{subequations}
	\label{eq:Sec-InfSupConstantFree-SubSec-InBandEigenmodes-FieldSolution}
	\begin{align}
	\label{eq:Sec-InfSupConstantFree-SubSec-InBandEigenmodes-FieldSolution1}
&\mathbf{\tilde E}(\omega) = i \omega \sum \limits_{\omega_n^2 \in \mathcal{B}_2} \frac{A_n}{1 - \frac{\omega^2}{\omega_n^2}} \mathbf{e}_n \text{, if } \omega^2 \neq \omega_n^2 \in \mathcal{B}_2,\\
	\label{eq:Sec-InfSupConstantFree-SubSec-InBandEigenmodes-FieldSolution2}
&\mathbf{\tilde E}(\omega) = i \omega \sum \limits_{ \omega_n^2 = \omega^2 } a_n \mathbf{e}_n + i \omega \sum \limits_{\omega_n^2 \in \mathcal{B}_2 \setminus \{\omega^2\}} \frac{A_n}{1 - \frac{\omega^2}{\omega_n^2}} \mathbf{e}_n \text{,} \\ &\text{ if } \omega^2 = \omega_n^2 \in \mathcal{B}_2 \text{,}
	\end{align}
\end{subequations}
and
\begin{subequations}
\label{eq:Sec-InfSupConstantFree-SubSec-InBandEigenmodes-ErrorSolution}
\begin{align}
	\label{eq:Sec-InfSupConstantFree-SubSec-InBandEigenmodes-ErrorSolution1}
&\boldsymbol{\tilde \epsilon}(\omega) = \mathbf{0} \text{, if } \omega^2 \neq \omega_n^2 \in \mathcal{B}_2,\\
	\label{eq:Sec-InfSupConstantFree-SubSec-InBandEigenmodes-ErrorSolution2}
&\boldsymbol{\tilde \epsilon}(\omega) = i \omega \sum \limits_{ \omega_n^2 = \omega^2 } a^{\boldsymbol{\epsilon}}_n \mathbf{e}_n \text{, if } \omega^2 = \omega_n^2 \in \mathcal{B}_2 \text{.}
	\end{align}
\end{subequations}
$\mathcal{B}_2$ stands for $[\omega^2_\text{min},\omega^2_\text{max}]$. As \eqref{eq:Sec-InfSupConstantFree-SubSec-InBandEigenmodes-ErrorSolution} shows, this is the worst-case scenario. Although the in-band eigenmode approximation basis in the eigenspace $\mathcal{V}_{\mathcal{B}}$ is the best basis to capture the fundamental dynamics of the electric field $\mathbf{E}(\omega)$ in the band of analysis $\mathcal{B}$ and so might be the case for the residual system, it turns out the opposite: the approximate state error $\boldsymbol{\tilde \epsilon}(\omega)$ is, apart from the in-band resonances, identically zero throughout the whole electromagnetic spectrum, $\forall \omega \in \mathbb{R}$ (see \eqref{eq:Sec-InfSupConstantFree-SubSec-InBandEigenmodes-ErrorSolution}). Unfortunately, the actual error is not zero but, on the contrary, we get zero error in \eqref{eq:Sec-InfSupConstantFree-SubSec-InBandEigenmodes-ErrorSolution}. This is how Galerkin approximation works, i.e., as far as the projection space is concerned, no error can be identified within this space since all the error, which is not identically zero, is orthogonal to the testing space and, therefore, remains outside the projection space used in \eqref{eq:Sec-InfSupConstantFree-SubSec-InBandEigenmodes-EigenBasis}. As a result, we get zero error in \eqref{eq:Sec-InfSupConstantFree-SubSec-InBandEigenmodes-ErrorSolution}, apart from the in-band resonances. This is the rationale behind the use of different projection spaces for the reduced problems \eqref{eq:Sec-InfSupConstantFree-GalerkinProjection} in \cite{feng2020InfSupConstantFree}, despite the fact that two similar problems have to be solved independently at the same time. Using two different projection spaces for the reduced systems increases the chance to identify the true state error in the approximation. In this work, we aim to achieve this goal while keeping the computational burden even much lower.
\subsection{Enhanced Reduced-Basis Space}
\label{Sec-InfSupConstantFree-SubSec-EnhancedReducedBasisSpace}
In-band eigenmodes are not enough to ensure a good approximation to the electric field $\mathbf{E}(\omega)$ in $\mathcal{B}$. In particular, only resonant phenomena are strictly captured by the in-band eigenmode basis in the eigenspace $\mathcal{V}_{\mathcal{B}}$, while other electromagnetic phenomena, such as direct source to load couplings, are missing. This has already been shown in Section \ref{Sec-InfSupConstantFree-SubSec-InBandEigenmodes}. As a result, the reduced-basis space has to be enriched by snapshots of the electric field $\mathbf{E}(w_k)$ in the frequency band $\mathcal{B}$, namely, 
\begin{subequations}
	\label{eq:Sec-InfSupConstantFree-SubSec-EnhancedReducedBasisSpace-RBS}
	\begin{align}
	\label{eq:Sec-InfSupConstantFree-SubSec-EnhancedReducedBasisSpace-RBS1}
	\mathcal{H}_m &= \mathcal{V}_\mathcal{B} + \text{span} \{\mathbf{E}(w_k)~|~w_k \in~\mathcal{B}\},\\
	\label{eq:Sec-InfSupConstantFree-SubSec-EnhancedReducedBasisSpace-RBS2}
	\mathcal{V}_\mathcal{B} &= \text{span}\{\mathbf{e}_n~|~\omega_n \in~\mathcal{B}\},
	\end{align}
\end{subequations}
giving rise to an enhanced reduced-basis space which ensures convergence to the electric field $\mathbf{E}(\omega)$ in $\mathcal{B}$ within a fast setting \cite{delaRubia2018CRBM}. In our situation, one question still remains: what is the reduced-basis space for the residual problem $\mathcal{H}^{\boldsymbol{\epsilon}}_m$ that should be used to get an approximation $\boldsymbol{\tilde \epsilon}(\omega)$ to the state error while keeping a low computational effort? We have already shown in Section \ref{Sec-InfSupConstantFree-SubSec-InBandEigenmodes} that using the same reduced-basis space considered for the original problem, i.e., $\mathcal{H}^{\boldsymbol{\epsilon}}_m = \mathcal{H}_m$, while keeping the computational burden low, is not a good choice (see \eqref{eq:Sec-InfSupConstantFree-SubSec-InBandEigenmodes-EigenBasis}--\eqref{eq:Sec-InfSupConstantFree-SubSec-InBandEigenmodes-ErrorSolution}). Also, carrying out a totally different reduced-basis space strategy for both reduced systems, i.e., the original and residual systems, yields good approximation results but the computational cost substantially increases \cite{feng2020InfSupConstantFree}. We then follow a compromise criterion: we allow the reduced-basis space for the residual problem $\mathcal{H}^{\boldsymbol{\epsilon}}_m$ to be different from the reduced-basis space for the original problem $\mathcal{H}_m$, but just for only one additional basis vector. As a result, we keep the computational cost low enough since only one additional solution needs to be carried out in the whole process. This cheap approach may result in underestimation of $\boldsymbol{\tilde \epsilon}(\omega)$ for the actual value of the state error $\boldsymbol{\epsilon}(\omega)$. This is not critical, as will become clear later. However, what is indeed essential is to choose an additional basis vector that allows to monitor the state error by means of the reduced residual system while not interfering the greedy algorithm with some unwanted frequency modulation in the snapshot selection for the reduced original problem. In other words, we should be careful to prevent the residual system solution from biasing the greedy procedure to solve the reduced original problem. Otherwise, the snapshots selected in the greedy process will end up with a bad choice for approximation purposes.

We propose to add the stationary electric field $\mathbf{F}_0$ as the only additional basis vector into $\mathcal{H}^{\boldsymbol{\epsilon}}_m$, as $\mathbf{F}_0$ is missing in $\mathcal{H}_m$, so it contributes to the error of the ROM. It is therefore reasonable to add $\mathbf{F}_0$ to the reduced basis space $\mathcal{H}^{\boldsymbol{\epsilon}}_m$ for the approximate error $\boldsymbol{\tilde \epsilon}(\omega)$. Furthermore, not only is it orthogonal to the eigenspace $\mathcal{V}_{\mathcal{B}}$ but also it has an almost-flat smooth influence on the fields $\mathbf{E}(\omega)$ throughout the frequency band of interest $\mathcal{B}$, due to its $1/\omega$ frequency behavior, cf. \eqref{eq:Sec-ProblemStatement-SolutionMaxwellFrequencyDependency}. As a consequence, its contribution to the error is also small, not spoiling the desired property for $\boldsymbol{\tilde \epsilon}(\omega)$, which should allow us to identify the missing information in the reduced original system. Thus, we can ensure that we are not modulating the adaptive sampling in the greedy algorithm. As a result, the basis vector $\mathbf{F}_0$ added to the reduced-basis space for the original problem $\mathcal{H}_m$ is a suitable candidate, since it has the desired behavior as testing space $\mathcal{H}^{\boldsymbol{\epsilon}}_m$ for the residual problem. Putting everything together, we use the following spaces:
\begin{subequations}
	\label{eq:Sec-InfSupConstantFree-SubSec-EnhancedReducedBasisSpace-ResidualRBS}
	\begin{align}
	\label{eq:Sec-InfSupConstantFree-SubSec-EnhancedReducedBasisSpace-ResidualRBS2}
	\mathcal{H}_m &= \mathcal{V}_\mathcal{B} + \text{span} \{\mathbf{E}(w_k)~|~w_k \in~\mathcal{B}\},\\
	\label{eq:Sec-InfSupConstantFree-SubSec-EnhancedReducedBasisSpace-ResidualRBS3}
	\mathcal{V}_\mathcal{B} &= \text{span}\{\mathbf{e}_n~|~\omega_n \in~\mathcal{B}\},\\
	\label{eq:Sec-InfSupConstantFree-SubSec-EnhancedReducedBasisSpace-ResidualRBS1}
	\mathcal{H}^{\boldsymbol{\epsilon}}_m &= \mathcal{H}_m + \text{span} \{ \mathbf{F}_0 \} \text{.}
	\end{align}
\end{subequations}
Unfortunately, the actual value of the approximate state error $\boldsymbol{\tilde \epsilon}(\omega)$ in this cheap statics-based reduced residual problem solution, while accurately identifying the missing information in the original problem, is \emph{not} a reliable indicator of the true state error $\boldsymbol{\epsilon}(\omega)$ due to its implicit underestimation, and can only be considered as a \emph{rough} estimator as compared to the one in Algorithm \ref{alg:Sec-InfSupConstantFree-ROMErrorEstimator}. As a result, a different strategy should be carried out to estimate the actual state error in the system.

We propose Algorithm \ref{alg:Sec-InfSupConstantFree-SubSec-EnhancedReducedBasisSpace-FastROMErrorEstimator} to carry out a fast \emph{a posteriori} state error estimation for reliable frequency sweep analysis in microwave devices. The reduced-basis space is adaptively built up by means of a greedy algorithm based on this state error estimation. Further, to overcome the lack of reliability of the rough error estimator $\boldsymbol{\tilde \epsilon}(\omega)$ in case it is used as the stopping criterion, a density stopping criterion \cite{delaRubia2018CRBM} is preferred in Step \ref{step:Sec-InfSupConstantFree-SubSec-EnhancedReducedBasisSpace-Alg-FastROMErrorEstimator-StoppingCriterion}. In contrast to the error estimator $\boldsymbol{\tilde \epsilon}(\omega^*)$ in Step \ref{step:Sec-InfSupConstantFree-Alg-ROMErrorEstimator-StoppingCriterion} of Algorithm \ref{alg:Sec-InfSupConstantFree-ROMErrorEstimator}, the true error at $\omega^*$ is computed in Step \ref{step:Sec-InfSupConstantFree-SubSec-EnhancedReducedBasisSpace-Alg-FastROMErrorEstimator-StoppingCriterion} of Algorithm \ref{alg:Sec-InfSupConstantFree-SubSec-EnhancedReducedBasisSpace-FastROMErrorEstimator} and is used as the stopping criterion. Note that the true error $\mathbf{E}(\omega^*)- \mathbf{\tilde E}(\omega^*)$ at $\omega^*$ is actually equivalent to the error $\mathbf{e}_\perp(\omega^*)$ defined in Algorithm \ref{alg:Sec-InfSupConstantFree-SubSec-EnhancedReducedBasisSpace-NewInformationNorm} and is readily available from Step \ref{step:Sec-InfSupConstantFree-SubSec-EnhancedReducedBasisSpace-Alg-FastROMErrorEstimator-NewSnapshot} in Algorithm \ref{alg:Sec-InfSupConstantFree-SubSec-EnhancedReducedBasisSpace-FastROMErrorEstimator}, where $\mathbf{E}(\omega^*)$ is orthogonalized against the existing basis vectors in $\mathcal{H}_m$ before being added to $\mathcal{H}_m$. In other words, the new error indicator in Step \ref{step:Sec-InfSupConstantFree-SubSec-EnhancedReducedBasisSpace-Alg-FastROMErrorEstimator-StoppingCriterion} of Algorithm \ref{alg:Sec-InfSupConstantFree-SubSec-EnhancedReducedBasisSpace-FastROMErrorEstimator} is based on how much new information the new greedy field sample $\mathbf{E}(\omega^*)$ is adding to the reduced-basis space. Eventually, the procedure stops whenever there is nothing new to add, within a tolerance denoted by \texttt{tol}, and the reduced-basis space can be considered as a dense enough approximation space which accurately describes the field solution $\mathbf{E}(\omega)$ throughout the band of interest $\mathcal{B}$. 

In summary, we propose Algorithm \ref{alg:Sec-InfSupConstantFree-SubSec-EnhancedReducedBasisSpace-FastROMErrorEstimator} which makes use of two error indicators: one is the rough error estimator $\boldsymbol{\tilde \epsilon}(\omega)$ employed in Step \ref{step:Sec-InfSupConstantFree-SubSec-EnhancedReducedBasisSpace-Alg-FastROMErrorEstimator-RoughErrorEstimator} based on a fixed $\mathcal{H}^{\mathbf{r}}_m$; the other is the indicator used in Step \ref{step:Sec-InfSupConstantFree-SubSec-EnhancedReducedBasisSpace-Alg-FastROMErrorEstimator-StoppingCriterion}, which is the true error at $\omega^*$ selected by the rough error estimator. It is then essential that these adaptive greedy field snapshots are properly selected in the band of analysis. Otherwise the procedure stops with no control on the actual error in the ROM. Fortunately, the rough error estimator $\boldsymbol{\tilde \epsilon}(\omega)$ is capable of providing what is indeed the missing information in the system. That means it catches the trend of the error change during the greedy algorithm, although it may underestimate the actual error. Therefore, we expect good performance in the proposed methodology. Finally, we present some remarks about Algorithm \ref{alg:Sec-InfSupConstantFree-SubSec-EnhancedReducedBasisSpace-FastROMErrorEstimator}.
\begin{algorithm}[tbp]
	\caption{Fast \emph{a posteriori} state error estimation for reliable frequency-parameter MOR.}
	\label{alg:Sec-InfSupConstantFree-SubSec-EnhancedReducedBasisSpace-FastROMErrorEstimator}
	\begin{algorithmic}[1]
		\REQUIRE Frequency band of interest $\mathcal{B}:=[\omega_\text{min}, \omega_\text{max} ]$, tolerance $\texttt{tol} > 0$ as the acceptable state error.
		\ENSURE $\mathcal{H}_m$ to ensure $\texttt{tol}$ state error in \eqref{eq:Sec-InfSupConstantFree-GalerkinProjection}.
		\STATE Solve for the in-band eigenmodes $\mathbf{e}_n$ in \eqref{eq:Sec-ProblemStatement-Subsec-StandarAPosterioriError-ParametricWeakForm}. Form the in-band eigenspace $\mathcal{V}_{\mathcal{B}}$: $\mathcal{V}_\mathcal{B} = \text{span}\{\mathbf{e}_n~|~\omega_n \in~\mathcal{B}\}$.				
		\STATE Solve for the stationary solution $\mathbf{F}_0$ in \eqref{eq:Sec-ProblemStatement-Subsec-StandarAPosterioriError-ParametricWeakForm}.					
		\STATE Initialize $\mathcal{H}_m = \mathcal{V}_\mathcal{B}$, $\mathcal{H}^{\mathbf{r}}_m = \text{span}\{ \mathbf{F}_0 \}$, $\xi = \texttt{tol} + 1$. Choose a sample $\omega^*$ randomly taken from the end points in $\mathcal{B}$.				
		\label{step:Sec-InfSupConstantFree-SubSec-EnhancedReducedBasisSpace-Alg-FastROMErrorEstimator-EndPoints}
		\STATE Solve for $\mathbf{E}(\omega^*)$ in \eqref{eq:Sec-ProblemStatement-Subsec-StandarAPosterioriError-ParametricWeakForm} and update $\mathcal{H}_m$: $\mathcal{H}_m = \mathcal{H}_m + \text{span}\{\mathbf{E}(\omega^*)\} $.
		\WHILE{$\xi > \texttt{tol}$}
		\STATE Compute the approximate field $\mathbf{\tilde E}(\omega)$ in \eqref{eq:Sec-InfSupConstantFree-GalerkinProjectionField}.
		\STATE Form $\mathcal{H}^{\boldsymbol{\epsilon}}_m$: $\mathcal{H}^{\boldsymbol{\epsilon}}_m = \mathcal{H}_m + \mathcal{H}^{\mathbf{r}}_m$.
		\STATE Compute the state error estimation $\boldsymbol{\tilde \epsilon}(\omega)$ in \eqref{eq:Sec-InfSupConstantFree-GalerkinProjectionError}.
		\label{step:Sec-InfSupConstantFree-SubSec-EnhancedReducedBasisSpace-Alg-FastROMErrorEstimator-RoughErrorEstimator}				
		\STATE Choose the next sample $\omega^*$ from $\mathcal{B}$ as $$\omega^*=\arg\max\limits_{\omega\in \mathcal{B}}\|\boldsymbol{\tilde \epsilon}(\omega)\|_{\mathcal{H}}.$$
		\STATE Compute $\mathbf{E}(\omega^*)$ in \eqref{eq:Sec-ProblemStatement-Subsec-StandarAPosterioriError-ParametricWeakForm} and update $\mathcal{H}_m$: $\mathcal{H}_m = \mathcal{H}_m + \text{span}\{\mathbf{E}(\omega^*)\} $.
		\label{step:Sec-InfSupConstantFree-SubSec-EnhancedReducedBasisSpace-Alg-FastROMErrorEstimator-NewSnapshot}
		\STATE $\xi = \|\mathbf{e}_\perp(\omega^*)\|_{\mathcal{H}} = \|\mathbf{E}(\omega^*)\|_{\nu}.$
		\label{step:Sec-InfSupConstantFree-SubSec-EnhancedReducedBasisSpace-Alg-FastROMErrorEstimator-StoppingCriterion}
		\ENDWHILE
		\STATE Use $\mathcal{H}_m$ to solve \eqref{eq:Sec-InfSupConstantFree-GalerkinProjectionField}.
	\end{algorithmic}
\end{algorithm}
\begin{algorithm}[tbp]
	\caption{Evaluation of the new information added to the reduced-basis space $\mathcal{H}_m$ by the field sample $\mathbf{E}$ \cite{delaRubia2018CRBM}.}
	\label{alg:Sec-InfSupConstantFree-SubSec-EnhancedReducedBasisSpace-NewInformationNorm}
	\begin{algorithmic}[1]
		\REQUIRE Electric field $\mathbf{E}$ and reduced-basis space $\mathcal{H}_m$.
		\ENSURE Norm $\|\mathbf{E}\|_{\nu}$ of the field  $\mathbf{E}$ indicating the new information added to the reduced-basis space $\mathcal{H}_m$.
		\STATE Normalize $\mathbf{E}$: $\mathbf{e} = \mathbf{E}/ \|\mathbf{E}\|_{\mathcal{H}}$.
		\STATE Project $\mathbf{e}$ onto $\mathcal{H}_m$. Decompose $\mathbf{e}$ into its projection and its orthogonal complement onto $\mathcal{H}_m$: $\mathbf{e} = \mathbf{e}_\parallel + \mathbf{e}_\perp $.
		\STATE Set $\|\mathbf{E}\|_{\nu} = \|\mathbf{e} - \mathbf{e}_\parallel\|_{\mathcal{H}} = \|\mathbf{e}_\perp\|_{\mathcal{H}}$.
	\end{algorithmic}	
\end{algorithm}
\begin{remark}
	Once the in-band eigenmode basis is completed, a randomly chosen sample between the end points in the band of analysis $\mathcal{B}$ is taken to enrich the eigenbasis in Step \ref{step:Sec-InfSupConstantFree-SubSec-EnhancedReducedBasisSpace-Alg-FastROMErrorEstimator-EndPoints}. At this point, neither the residual norm nor the state error norm can provide an answer to what is the best sample to choose, since both estimators show a constant behavior throughout the whole spectrum $\forall \omega \in \mathbb{R}$. We refer to \cite{delaRubia2018CRBM} for the illustration of the residual norm behavior when eigenbasis is used. As a result, based on the new information arguments added to the eigenbasis, any of the end points in the frequency band of interest should be sampled (see \eqref{eq:Sec-ProblemStatement-SolutionMaxwellFrequencyDependency}). It should be noted that, in \eqref{eq:Sec-ProblemStatement-SolutionMaxwellFrequencyDependency}, the further we are from the in-band eigenmodes, the more linearly independent new information is found, until a new out of band eigenmode eventually shows up. This is the main reason for sampling at the end points of the frequency band of interest $\mathcal{B}$.
\end{remark}
%
%
\begin{remark}
	Contrary to what is done in \cite{delaRubia2018CRBM}, the residual norm is not used to guide the greedy algorithm at any step. Thus, better results are expected. Residual information is problematic since its behavior suffers from the in-band eigenresonances. The residual does not vanish at and nearby resonances, which can definitely mislead the greedy algorithm. This will become apparent in Section \ref{Sec-NumericalResults} throughout the numerical examples.
\end{remark}
\section{Numerical Results}
\label{Sec-NumericalResults}
In this section, we apply the proposed \emph{a posteriori} state error estimation for reliable fast frequency sweeps of different challenging microwave circuits, namely, a quad-mode dielectric resonator filter, an inline filter with transmission zeros generated by frequency-dependent couplings, an inline dielectric resonator filter and a combline diplexer. The capabilities and reliability of the proposed procedure are demonstrated via these examples. The in-house C++ code for finite element method (FEM) simulations uses a second-order first family of N\'ed\'elec's elements \cite{Ned80, Ing06}, on meshes provided by \texttt{Gmsh} \cite{GeuR09}. All computations were carried out on a workstation with 3.00-GHz Intel Xeon E5-2687W v4 processor and 256-GB RAM. 

In our experiments, we define the true error ($\epsilon_{\text{true}}$) of the ROM as the maximal error over the whole frequency band $\mathcal{B}$ using the indicator $\mathbf{e}_\perp$ in Algorithm \ref{alg:Sec-InfSupConstantFree-SubSec-EnhancedReducedBasisSpace-NewInformationNorm}, namely,
\begin{equation}
\label{eq:Sec-NumericalResults-TrueError}
\epsilon_{\text{true}} = \max \limits_{ \omega \in \mathcal{B} } \| \mathbf{e}_\perp(\omega) \|_{\mathcal{H}} \text{.}
\end{equation}
Note that the novelty of the proposed Algorithm \ref{alg:Sec-InfSupConstantFree-SubSec-EnhancedReducedBasisSpace-FastROMErrorEstimator} is twofold: 1) a rough error estimator $\boldsymbol{\tilde \epsilon}(\omega)$ to save computational costs as compared with the one in Algorithm \ref{alg:Sec-InfSupConstantFree-ROMErrorEstimator}; 2) the indicator $\| \mathbf{e}_\perp(\omega) \|_{\mathcal{H}}$ to improve the reliability of the rough estimator. To show that the proposed Algorithm \ref{alg:Sec-InfSupConstantFree-SubSec-EnhancedReducedBasisSpace-FastROMErrorEstimator} is more efficient,  we compare it with Algorithm \ref{alg:Sec-InfSupConstantFree-ROMErrorEstimator}, as well as with the greedy algorithm using the residual norm-based error estimator. To this end, we define the indicators $\| \mathbf{e}_\perp(\omega) \|_{\mathcal{H}}$ based on the above three different estimators, respectively, i.e., 
\begin{equation}
\label{eq:Sec-NumericalResults-StateErrorEstimator}
\epsilon_{\text{state}} = \|\mathbf{e}_\perp(  \arg\max\limits_{\omega\in \mathcal{B}}\|\boldsymbol{\tilde \epsilon}(\omega)\|_{\mathcal{H}} )\|_{\mathcal{H}} \text{,}
\end{equation}
where $\boldsymbol{\tilde \epsilon}(\omega)$ refers to either the one in Algorithm \ref{alg:Sec-InfSupConstantFree-ROMErrorEstimator} or the rough estimator in Algorithm \ref{alg:Sec-InfSupConstantFree-SubSec-EnhancedReducedBasisSpace-FastROMErrorEstimator}. The indicator based on the residual-norm is defined as
\begin{equation}
\label{eq:Sec-NumericalResults-ResidualBasedErrorEstimator}
\epsilon_{\text{res}} = \|\mathbf{e}_\perp( \arg\max\limits_{\omega \in \mathcal{B}} \|r(\mathbf{\tilde E}(\omega), \cdot; \omega)\|_{\mathcal{H}^{\prime}} )\|_{\mathcal{H}} \text{.}
\end{equation}
For a fair comparison, we compare the results of Algorithm \ref{alg:Sec-InfSupConstantFree-SubSec-EnhancedReducedBasisSpace-FastROMErrorEstimator} with Algorithm \ref{alg:Sec-InfSupConstantFree-ROMErrorEstimator}, where $\xi$ in Step \ref{step:Sec-InfSupConstantFree-Alg-ROMErrorEstimator-StoppingCriterion} is replaced by the indicator $\epsilon_{\text{state}}$ in \eqref{eq:Sec-NumericalResults-StateErrorEstimator}. We also compare Algorithm \ref{alg:Sec-InfSupConstantFree-SubSec-EnhancedReducedBasisSpace-FastROMErrorEstimator} with the residual-norm based greedy algorithm, where the indicator in \eqref{eq:Sec-NumericalResults-ResidualBasedErrorEstimator} is used as the stopping criterion. The comparison between Algorithm \ref{alg:Sec-InfSupConstantFree-ROMErrorEstimator} and Algorithm \ref{alg:Sec-InfSupConstantFree-SubSec-EnhancedReducedBasisSpace-FastROMErrorEstimator} is carried out for the first two examples. The indicator $\mathbf{e}_\perp$ plays a role of paramount importance in providing a fair comparison among the different strategies. It should be pointed out that the true error $\epsilon_{\text{true}}$ in \eqref{eq:Sec-NumericalResults-TrueError} implies the computation of the field solution $\mathbf{E}(\omega)$ by means of time-consuming FEM simulations throughout the whole frequency band of interest $\mathcal{B}$. This can \emph{only} be carried out for academic purposes.
Finally, as a figure of merit, we use the metric of effectivity to gauge how close the estimated error is to the true error:
\mbox{$	\texttt{eff} := \frac{\epsilon}{\epsilon_{\text{true}}}$,}
$\epsilon$ refers to either $\epsilon_{\text{state}}$ in \eqref{eq:Sec-NumericalResults-StateErrorEstimator} or $\epsilon_{\text{res}}$ in \eqref{eq:Sec-NumericalResults-ResidualBasedErrorEstimator}. For each of the four microwave circuits considered, we evaluate the performance of the error estimators using the indicators defined in \eqref{eq:Sec-NumericalResults-StateErrorEstimator} and \eqref{eq:Sec-NumericalResults-ResidualBasedErrorEstimator}. A tolerance threshold of $\texttt{tol} = 2\cdot10^{-7}$ is used throughout all the numerical examples, for all greedy algorithms corresponding to their respective error estimators.
%
\subsection{Quad-Mode Dielectric Resonator Filter}
\label{Sec-NumericalResults-Subsec-QuadModeFilter}
\begin{figure}[tbp]
	\centering
	\includegraphics[width=0.7\linewidth]{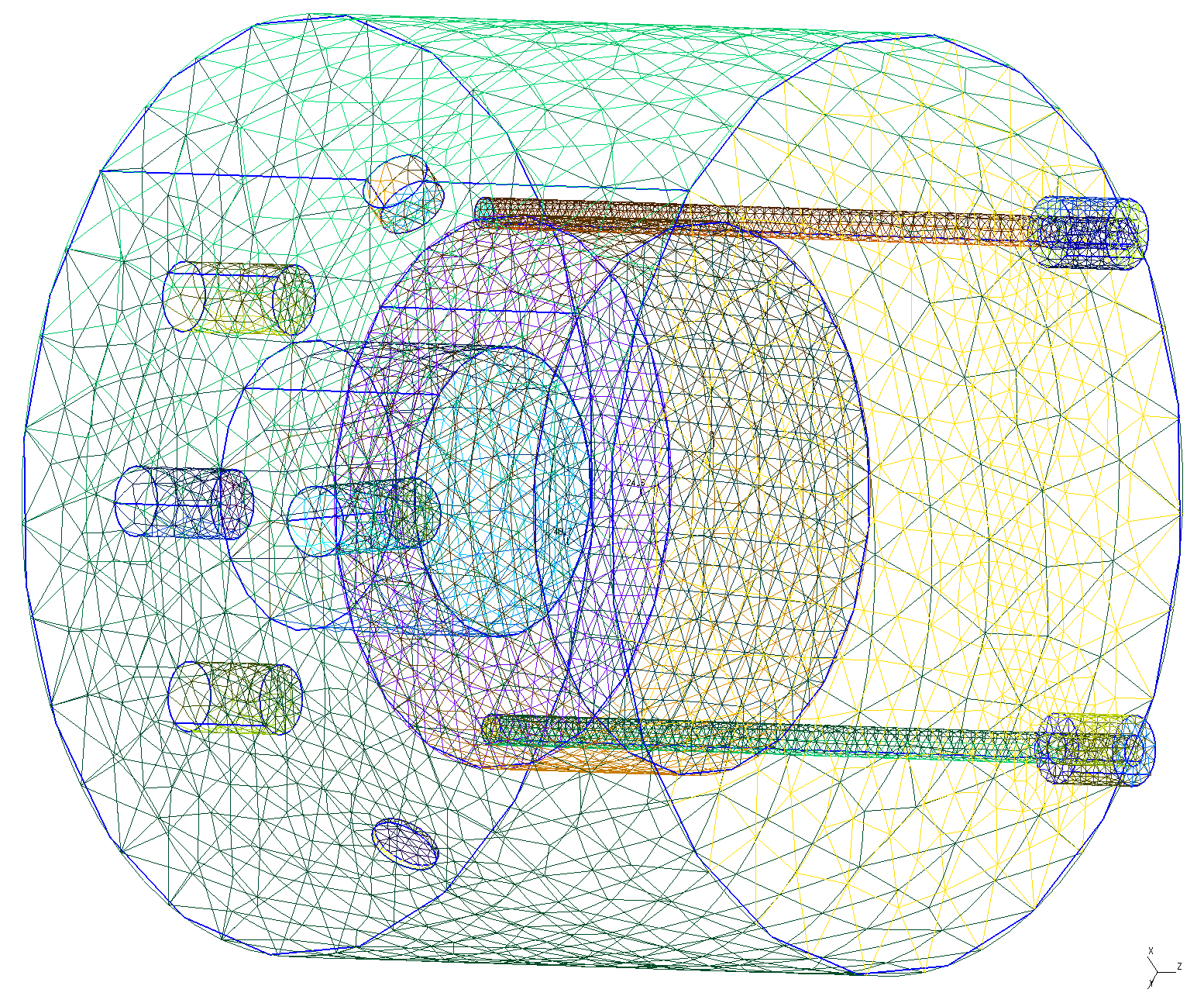}
	\caption{Quad-mode dielectric resonator filter proposed in \cite{memarian2009quad}.}
	\label{fig:Sec-NumericalResults-Subsec-QuadModeFilter-FilterGeometry}
\end{figure}

\begin{figure}[tbp]
	\centering
	\input{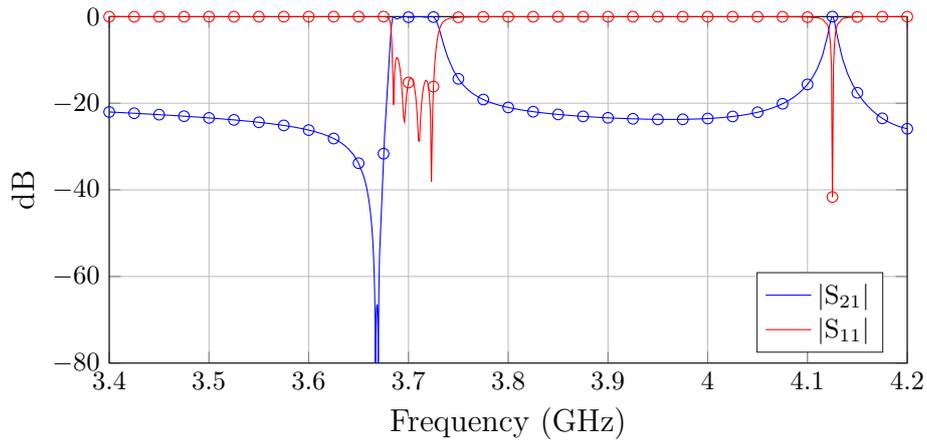}
	\caption{Quad-mode dielectric resonator filter scattering parameter response comparison with the proposed approach. (--) RBM. ($\circ$) FEM.}
	\label{fig:Sec-NumericalResults-Subsec-QuadModeFilter-FilterResponse}
\end{figure}

A quad-mode dielectric resonator filter in a single cylindrical cavity is proposed in \cite{memarian2009quad}. Fig. \ref{fig:Sec-NumericalResults-Subsec-QuadModeFilter-FilterGeometry} shows the geometry of the filter as well as the mesh used for its analysis. These structures are extremely attractive since multiple resonant modes show up in a single cylindrical cavity due to the dielectric resonator. At the same time, they are difficult to tune since all dominant modes are coupled with each other, requiring multiple full-wave electromagnetic analyses to carry out a good electrical design. Six tuning screws are included in this filter. It is then of paramount importance to accurately predict the electromagnetic behavior in a frequency band in an efficient way. The filter detailed in Fig. \ref{fig:Sec-NumericalResults-Subsec-QuadModeFilter-FilterGeometry} is a quad-mode filter, where a four-pole passband filtering response, shown in Fig. \ref{fig:Sec-NumericalResults-Subsec-QuadModeFilter-FilterResponse}, is obtained. However, at the same time, there are additional resonant modes in the band of analysis, $\mathcal{B}:=[3.4, 4.2]$ GHz, giving rise to a more complicated response including direct source to load coupling, affecting the position of the two transmission zeros, rather than a typical four-pole frequency response. As a result, a reliable ROM for fast frequency sweep analysis is essential.

An FEM system with 245,778 degrees of freedom is used to solve for the electric field. Following Algorithm \ref{alg:Sec-InfSupConstantFree-SubSec-EnhancedReducedBasisSpace-FastROMErrorEstimator}, a ROM is obtained by means of the Reduced-Basis Method (RBM) giving rise to a reduced system of dimension 14 to compute the frequency response detailed in Fig. \ref{fig:Sec-NumericalResults-Subsec-QuadModeFilter-FilterResponse}. It is clear that there is not even a competition ($14 \ll 245,778$) between the computational effort carried out by RBM and the one that would have been required by FEM to get the same frequency response along the same frequency samples, nevermind subsampling is considered. Solving a reduced system of dimension 14 many times is totally effortless while solving many times an FEM system of dimension 245,778 is rather time-consuming. Good agreement is found between the FEM and RBM results in Fig. \ref{fig:Sec-NumericalResults-Subsec-QuadModeFilter-FilterResponse}.

Next, a comparison of the different MOR methodologies is carried out. Table \ref{tab:Sec-NumericalResults-Subsec-QuadModeFilter-GreedySamples} not only details the frequency samples adaptively chosen by each greedy algorithm but also the estimated and true errors, at each iteration in the MOR process. Fig. \ref{fig:Sec-NumericalResults-Subsec-QuadModeFilter-Error} depicts the convergence behavior for the different methodologies and details the effectivity metric for the proposed \emph{a posteriori} state error estimator. Good behavior is observed in the proposed methodology in comparison to the true error. On the contrary, the residual norm-based error estimation results in a greedy algorithm that underestimates the error and prematurely stops the iterative procedure. The rationale behind this is shown in Fig. \ref{fig:Sec-NumericalResults-Subsec-QuadModeFilter-ErrorEstimator}. For a ROM of dimension 12 obtained by the proposed approach, both the residual norm $\|r(\mathbf{\tilde E}(\omega), \cdot; \omega)\|_{\mathcal{H}^{\prime}}$ and the approximate state error  $\|\boldsymbol{\tilde \epsilon}(\omega)\|_{\mathcal{H}}$ are plotted versus frequency. While the state error estimation has a smooth frequency behavior, the residual norm is contaminated by the in-band resonant modes, creating an eigenfrequency pollution. This makes the residual norm-based greedy algorithm mislead its sampling once again around the eigenresonances, notoriously deteriorating the new information added to the reduced-basis space. This is the key advantage of using state error estimation, where this unwanted behavior does not hold.

Finally, Table \ref{tab:Sec-NumericalResults-Subsec-QuadModeFilter-AlgorithmComparison} compares the performance of the proposed approach with Algorithm~\ref{alg:Sec-InfSupConstantFree-ROMErrorEstimator}, where $\xi$ in Step \ref{step:Sec-InfSupConstantFree-Alg-ROMErrorEstimator-StoppingCriterion} is replaced by the indicator in \eqref{eq:Sec-NumericalResults-StateErrorEstimator}. It should be pointed out that, even though the in-band eigenmodes are not imposed in Algorithm~\ref{alg:Sec-InfSupConstantFree-ROMErrorEstimator}, these in-band eigenmodes do show up in the greedy algorithm after the first few iterations in the procedure. This indicates Algorithm \ref{alg:Sec-InfSupConstantFree-ROMErrorEstimator} is accurately working out since the in-band eigenmodes have been shown to be the best choice from the theoretical point of view in Section \ref{Sec-ProblemStatement}. As far as computational complexity is concerned, Algorithm \ref{alg:Sec-InfSupConstantFree-ROMErrorEstimator} is more expensive than Algorithm \ref{alg:Sec-InfSupConstantFree-SubSec-EnhancedReducedBasisSpace-FastROMErrorEstimator}. As a matter of fact, the size of both ROMs when the same stopping criterion is used is $14$ for Algorithm \ref{alg:Sec-InfSupConstantFree-SubSec-EnhancedReducedBasisSpace-FastROMErrorEstimator} and $16$ for Algorithm \ref{alg:Sec-InfSupConstantFree-ROMErrorEstimator}, whereas the size of the residual ROMs defined in  \eqref{eq:Sec-InfSupConstantFree-GalerkinProjectionError} is $15$ and $30$ for Algorithms \ref{alg:Sec-InfSupConstantFree-SubSec-EnhancedReducedBasisSpace-FastROMErrorEstimator} and \ref{alg:Sec-InfSupConstantFree-ROMErrorEstimator}, respectively. This shows how many more efforts we need to carry out in Algorithm \ref{alg:Sec-InfSupConstantFree-ROMErrorEstimator} compared to Algorithm \ref{alg:Sec-InfSupConstantFree-SubSec-EnhancedReducedBasisSpace-FastROMErrorEstimator}. As a result, the proposed approach yields a fast \emph{a posteriori} state error estimation.

\begin{table}[tbp]
	\centering
	\begin{threeparttable}
		\renewcommand{\arraystretch}{1.3}
		\caption{Greedy Algorithm Frequency Samples for Different MOR Methodologies in the Quad-Mode Filter.}
		\label{tab:Sec-NumericalResults-Subsec-QuadModeFilter-GreedySamples}
		\centering
		\begin{tabular}{c|c|c|c|c|c}
			\hline
			GHz & $\epsilon_{\text{true}}$ & GHz & $\epsilon_{\text{state}}$ & GHz & $\epsilon_{\text{res}}$ \\ 
			\hline
			\hline
			3.6103     & 1. 				& 3.6103 	& 1. 				& 3.6103	& 1. \\
			3.6245     & 1. 				& 3.6245 	& 1. 				& 3.6245	& 1. \\
			3.6887     & 1. 				& 3.6887 	& 1. 				& 3.6887	& 1. \\
			3.7149     & 1. 				& 3.7149 	& 1. 				& 3.7149	& 1. \\
			3.9299     & 1. 				& 3.9299 	& 1. 				& 3.9299	& 1. \\
			4.0067     & 1. 				& 4.0067 	& 1. 				& 4.0067	& 1. \\
			4.1552     & 1. 				& 4.1552 	& 1. 				& 4.1552	& 1. \\
			4.2000     & $1.7\cdot10^{-1}$ 	& 3.4000	& $1.6\cdot10^{-1}$ & 3.4000	& $1.6\cdot10^{-1}$ \\	
			3.4000     & $7.1\cdot10^{-2}$ 	& 4.2000	& $8.0\cdot10^{-2}$ & 4.0067	& $2.8\cdot10^{-7}$ \\						
			3.7800     & $1.6\cdot10^{-3}$ 	& 3.8820 	& $8.3\cdot10^{-4}$ & --		& -- \\
			4.1000     & $9.4\cdot10^{-5}$ 	& 4.0940	& $6.7\cdot10^{-5}$ & --		& -- \\
			3.4900     & $6.8\cdot10^{-6}$ 	& 3.5930 	& $2.5\cdot10^{-6}$	& --		& -- \\
			4.1800     & $1.2\cdot10^{-6}$ 	& 4.1680 	& $5.8\cdot10^{-7}$	& --		& -- \\
			3.8600     & $1.9\cdot10^{-7}$	& 3.6930 	& $1.5\cdot10^{-7}$	& --		& -- \\
			\hline
		\end{tabular}
	\end{threeparttable}
\end{table}
\begin{table}[tbp]
	\centering
	\begin{threeparttable}
		\renewcommand{\arraystretch}{1.3}
		\caption{State Error Algorithm Comparison in the Quad-Mode Filter.}
		\label{tab:Sec-NumericalResults-Subsec-QuadModeFilter-AlgorithmComparison}
		\centering
		\begin{tabular}{c|c|c|c}
			\hline
			\multicolumn {2} {c|}{Algorithm \ref{alg:Sec-InfSupConstantFree-SubSec-EnhancedReducedBasisSpace-FastROMErrorEstimator}} & \multicolumn {2} {|c}{Algorithm \ref{alg:Sec-InfSupConstantFree-ROMErrorEstimator}} \\ 
			\hline
			\hline
			GHz & $\epsilon_{\text{state}}$ & GHz & $\epsilon_{\text{state}}$ \\ 
			\hline
			\hline
			3.6103 	& 1. 				& 3.4000	& $1.6\cdot10^{-3}$ \\
			3.6245 	& 1. 				& 3.6470	& $1.7\cdot10^{-3}$ \\
			3.6887 	& 1. 				& 4.0810	& $1.6$ \\
			3.7149 	& 1. 				& 3.9300	& $2.6\cdot10^{-3}$ \\
			3.9299 	& 1. 				& 3.6250	& $1.2\cdot10^{-3}$ \\
			4.0067 	& 1. 				& 3.7150	& $3.2\cdot10^{-4}$ \\
			4.1552 	& 1. 				& 4.1550	& $9.8\cdot10^{-7}$ \\
			3.4000	& $1.6\cdot10^{-1}$ & 4.0070	& $5.5\cdot10^{-8}$ \\
			4.2000	& $8.0\cdot10^{-2}$ & --		& -- \\
			3.8820 	& $8.3\cdot10^{-4}$ & --		& -- \\
			4.0940	& $6.7\cdot10^{-5}$ & --		& -- \\
			3.5930 	& $2.5\cdot10^{-6}$	& --		& -- \\
			4.1680 	& $5.8\cdot10^{-7}$	& --		& -- \\
			3.6930 	& $1.5\cdot10^{-7}$	& --		& -- \\
			\hline
			\hline
			ROM size 			& 14 	& ROM size			& 16 \\
			Residual ROM size 	& 15 	& Residual ROM size	& 30 \\
			\hline
		\end{tabular}
	\end{threeparttable}
\end{table}
%

%

\begin{figure}[tbp]
	\centering
	\subfloat[]{\label{fig:Sec-NumericalResults-Subsec-QuadModeFilter-ErrorConvergence}
%
%
\definecolor{mycolor1}{rgb}{0.00000,0.44700,0.74100}%
\begin{tikzpicture}

\begin{axis}[%
width=0.31\fwidth,
height=0.31\fheight,
at={(0\fwidth,0\fheight)},
scale only axis,
xmin=0,
xmax=14,
xtick = {2,4,6,8,10,12,14},
xlabel style={font=\color{white!15!black}, font = \large},
xlabel={Iterations},
ymode=log,
ymin=1e-7,
ymax=1e0,
yminorticks=true,
ytick = {1e0, 1e-2, 1e-4, 1e-6},
ylabel style={font=\color{white!15!black}, font = \large},
ylabel={Error},
axis background/.style={fill=white},
xmajorgrids,
ymajorgrids,
yminorgrids,
legend pos=south west,
legend style={at={(0.01,0.1)},anchor=south west, legend cell align=left, align=left, draw=white!15!black, font=\small}
]
\addplot [color=blue, thick]
  table[row sep=crcr]{%
1	1.0\\
2	1.0\\
3	1.0\\
4	1.0\\
5	1.0\\
6	1.0\\
7	1.0\\
8   0.157761528897257\\
9	0.000000287183739\\
};
\addlegendentry{$\epsilon_{\text{res}}$}

\addplot [color=red, thick]
table[row sep=crcr]{%
1	1.0\\
2	1.0\\
3	1.0\\
4	1.0\\
5	1.0\\
6	1.0\\
7	1.0\\
8   0.157761528897257\\
9   0.079797368377660\\
10  0.000832941774676\\
11  0.000067117136411\\
12  0.000002497506757\\
13  0.000000582492060\\
14  0.000000152190998\\
};
\addlegendentry{$\epsilon_{\text{state}}$}

\addplot [color=black, dashed, mark=diamond, mark options={solid, black}, thick]
  table[row sep=crcr]{%
1	1.0\\
2	1.0\\
3	1.0\\
4	1.0\\
5	1.0\\
6	1.0\\
7	1.0\\
8   0.176789451042759\\
9   0.070710678118655\\
10  0.001577003804688\\
11  0.000093651775210\\
12  0.000006791747934\\
13  0.000001230058942\\
14  0.000000194439671\\
};
\addlegendentry{$\epsilon_{\text{true}}$}

\end{axis}
\end{tikzpicture}%
} 
		\subfloat[]{\label{fig:Sec-NumericalResults-Subsec-QuadModeFilter-ErrorEffectivity}
%
%
\definecolor{mycolor1}{rgb}{0.00000,0.44700,0.74100}%
\begin{tikzpicture}

\begin{axis}[%
width=0.31\fwidth,
height=0.31\fheight,
at={(0\fwidth,0\fheight)},
scale only axis,
xmin=1,
xmax=14,
xtick = {2,4,6,8,10,12,14},
xlabel style={font=\color{white!15!black}, font = \large},
xlabel={Iterations},
yminorticks=true,
ylabel style={font=\color{white!15!black}, font = \large},
ylabel={\texttt{eff}},
axis background/.style={fill=white},
xmajorgrids,
ymajorgrids,
yminorgrids,
legend pos=south west,
legend style={at={(0.01,0.1)},anchor=south west, legend cell align=left, align=left, draw=white!15!black, font=\small}
]
\addplot [color=red, thick]
table[row sep=crcr]{%
1   1.000000000000000 \\
2   1.000000000000000 \\
3   1.000000000000000 \\
4   1.000000000000000 \\
5   1.000000000000000 \\
6   1.000000000000000 \\
7   1.000000000000000 \\
8   0.892369584082821 \\
9   1.128505206013690 \\
10  0.528179939832562 \\
11  0.716666996010966 \\
12  0.367726656081810 \\
13  0.473548088059499 \\
14  0.7827157761442622\\
};
\addlegendentry{$\epsilon_{\text{state}}$}

\end{axis}
\end{tikzpicture}%
}
	\caption{Quad-mode dielectric resonator filter error estimator results. (a) Convergence of the greedy algorithm. (b) Effectivity (\texttt{eff}).}
	\label{fig:Sec-NumericalResults-Subsec-QuadModeFilter-Error}
\end{figure}
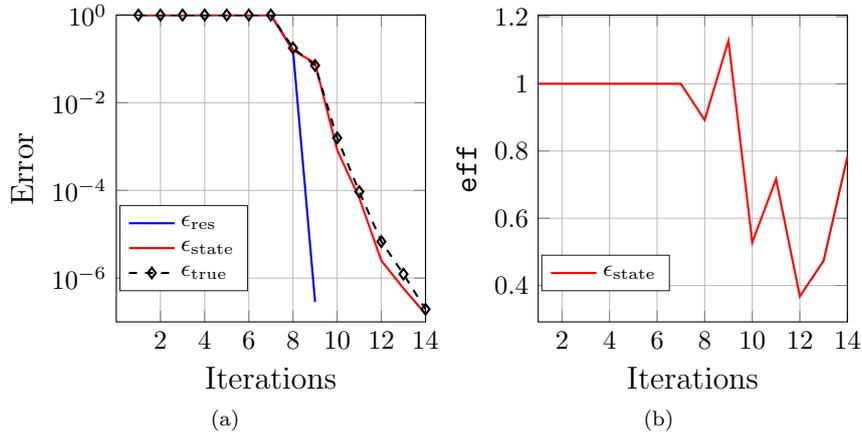	
\begin{figure}[tbp]
	\centering
	\input{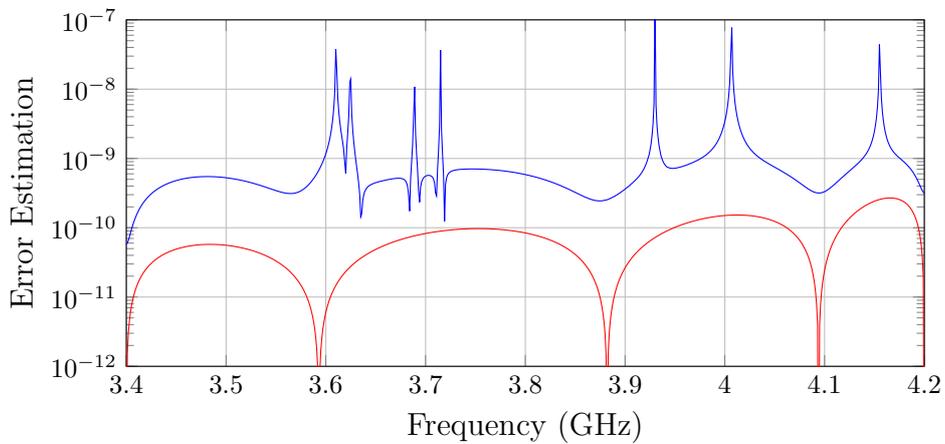}
	\caption{Quad-mode dielectric resonator filter error estimation frequency behavior for the ROM of order 12. ({\color{blue}{--}}) Residual error $\|r(\mathbf{\tilde E}(\omega), \cdot; \omega)\|_{\mathcal{H}^{\prime}}$. ({\color{red}{--}}) State error estimator $\|\boldsymbol{\tilde \epsilon}(\omega)\|_{\mathcal{H}}$.}
	\label{fig:Sec-NumericalResults-Subsec-QuadModeFilter-ErrorEstimator}
\end{figure}	
\subsection{Inline Filter with Frequency-Dependent Couplings}
\label{Sec-NumericalResults-Subsec-MacchiarellaFilter}
\begin{figure}[tbp]
	\centering
	\includegraphics[width=\linewidth]{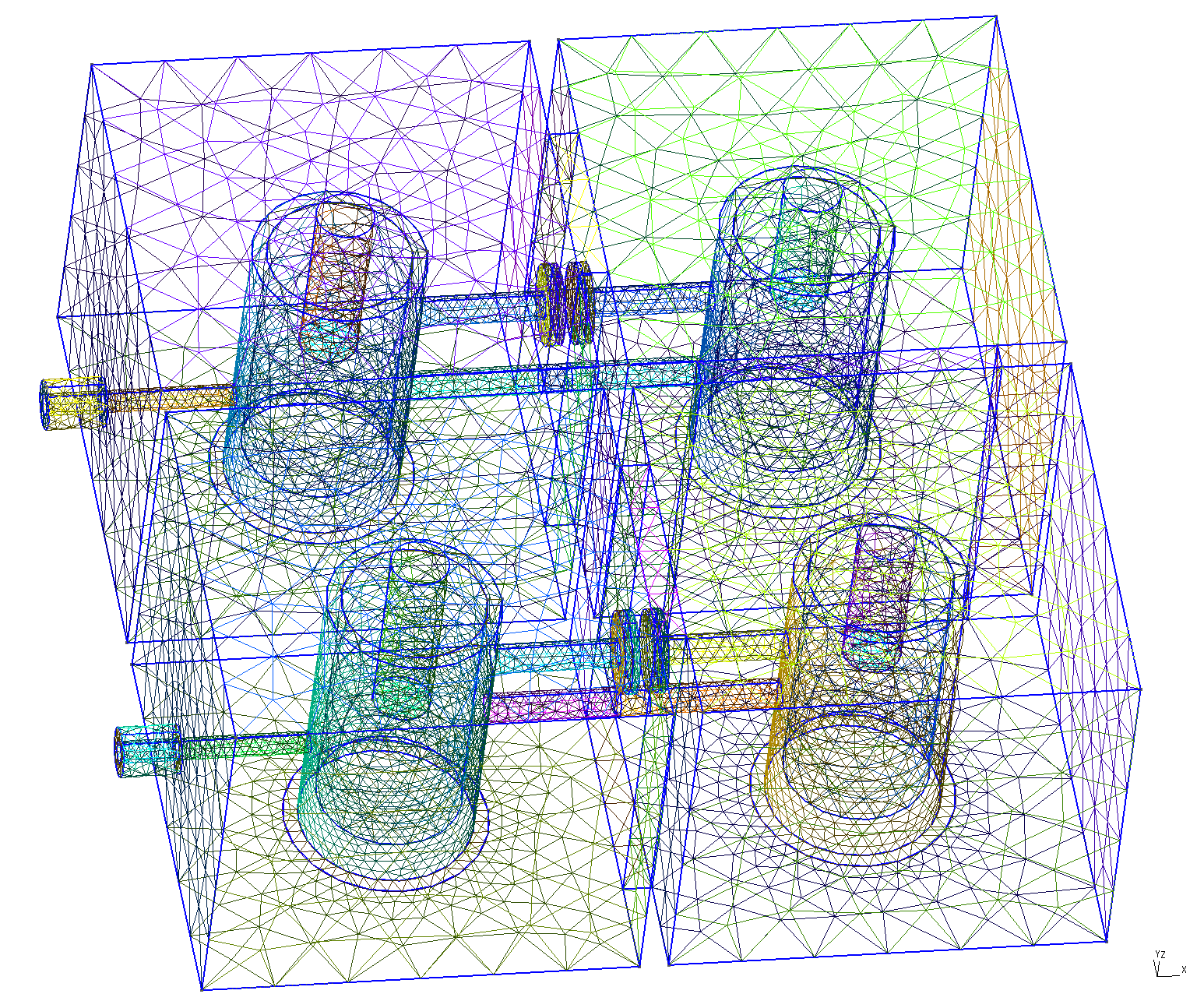}
	\caption{Inline filter with finite transmission zeros proposed in \cite{he2018direct}.}
	\label{fig:Sec-NumericalResults-Subsec-MacchiarellaFilter-FilterGeometry}
\end{figure}

The next example is a fourth-order inline combline filter designed in \cite{he2018direct}, where frequency dependent couplings are taken into account to provide finite transmission zeros, even within an inline coupling route structure. The filter geometry is depicted in Fig. \ref{fig:Sec-NumericalResults-Subsec-MacchiarellaFilter-FilterGeometry} and the mesh for the electromagnetic analysis is shown as well. Within four combline resonant cavities, two additional shunt inductors and capacitors are included, which indeed give rise to additional higher frequency resonances creating two finite transmission zeros near by the filter passband. As a result, better rejection properties are allowed in this inline structure. Frequency dependent coupling filter theory is continuously developing \cite{szydlowski2012coupled, tamiazzo2017synthesis} and a large number of full-wave analyses are required to design these kind of filtering responses. It is essential to carry out a reliable fast frequency sweep analysis to meet the target electrical response in the design optimization loop.

\begin{figure}[tbp]
	\centering
	\input{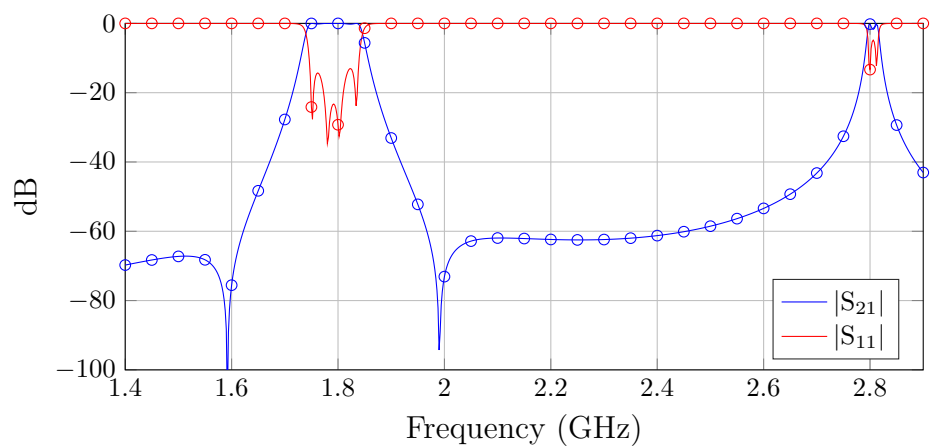}
	\caption{Inline filter with finite transmission zeros scattering parameter response comparison with the proposed approach. (--) RBM. ($\circ$) FEM.}
	\label{fig:Sec-NumericalResults-Subsec-MacchiarellaFilter-FilterResponse}
\end{figure}

The frequency band $\mathcal{B}:=[1.4, 2.9]$ GHz is taken into account for analysis. We solve for the electric field in the band of interest $\mathcal{B}$ by means of an FEM system with 105,690 degrees of freedom. Following the procedure proposed in this work, we obtain a ROM of dimension 13 to sweep the frequency response of this filter via RBM. A comparison between the filter response results obtained by FEM and RBM is shown in Fig. \ref{fig:Sec-NumericalResults-Subsec-MacchiarellaFilter-FilterResponse}. Reasonable agreement is achieved.

\begin{figure}[tbp]
	\centering
	\subfloat[]{\label{fig:Sec-NumericalResults-Subsec-MacchiarellaFilter-ErrorConvergence}
%
%
\definecolor{mycolor1}{rgb}{0.00000,0.44700,0.74100}%
\begin{tikzpicture}

\begin{axis}[%
width=0.3\fwidth,
height=0.3\fheight,
at={(0\fwidth,0\fheight)},
scale only axis,
xmin=1,
xmax=13,
xtick = {2,4,6,8,10,12},
xlabel style={font=\color{white!15!black}, font = \large},
xlabel={Iterations},
ymode=log,
ymax=1e0,
yminorticks=true,
ytick = {1e0, 1e-2, 1e-4, 1e-6, 1e-8},
ylabel style={font=\color{white!15!black}, font = \large},
ylabel={Error},
axis background/.style={fill=white},
xmajorgrids,
ymajorgrids,
yminorgrids,
legend pos=south west,
legend style={legend cell align=left, align=left, draw=white!15!black}
]
\addplot [color=blue, thick]
  table[row sep=crcr]{%
1	1.0\\
2	1.0\\
3	1.0\\
4	1.0\\
5	1.0\\
6	1.0\\
7   0.673492390454413 \\
8   0.153504071607238 \\
9  	0.000000084086681 \\
};
\addlegendentry{$\epsilon_{\text{res}}$}

\addplot [color=red, thick]
table[row sep=crcr]{%
1	1.0\\
2	1.0\\
3	1.0\\
4	1.0\\
5	1.0\\
6	1.0\\
7   0.673492390454413\\
8   0.153504723054374\\
9  	0.002775305388601\\
10  0.000348387715053\\
11  0.000007601124917\\
12  0.000000009419692\\
};
\addlegendentry{$\epsilon_{\text{state}}$}


\addplot [color=black, dashed, mark=diamond, mark options={solid, black}, thick]
  table[row sep=crcr]{%
1	1.0\\
2	1.0\\
3	1.0\\
4	1.0\\
5	1.0\\
6	1.0\\
7  	0.880055680056666\\
8  	0.077511289500304\\
9  	0.004683877667062\\
10 	0.000294652371448\\
11 	0.000009830368254\\
12 	0.000000490649264\\
13 	0.000000171538246\\
};
\addlegendentry{$\epsilon_{\text{true}}$}

\end{axis}
\end{tikzpicture}%
} 
	\subfloat[]{\label{fig:Sec-NumericalResults-Subsec-MacchiarellaFilter-ErrorEffectivity}
%
%
\definecolor{mycolor1}{rgb}{0.00000,0.44700,0.74100}%
\begin{tikzpicture}

\begin{axis}[%
width=0.3\fwidth,
height=0.3\fheight,
at={(0\fwidth,0\fheight)},
scale only axis,
xmin=1,
xmax=13,
xtick = {2,4,6,8,10,12},
xlabel style={font=\color{white!15!black}, font = \large},
xlabel={Iterations},
yminorticks=true,
ylabel style={font=\color{white!15!black}, font = \large},
ylabel={\texttt{eff}},
axis background/.style={fill=white},
xmajorgrids,
ymajorgrids,
yminorgrids,
legend pos=south west,
legend style={legend cell align=left, align=left, draw=white!15!black}
]
\addplot [color=red, thick]
table[row sep=crcr]{%
	1	   1.000000000000000 \\
	2	   1.000000000000000 \\
	3	   1.000000000000000 \\
	4	   1.000000000000000 \\
	5	   1.000000000000000 \\
	6	   1.000000000000000 \\
	7	   0.765283840235027 \\
	8	   1.980417614569178 \\
	9	   0.592523030248648 \\
	10	   1.182368610647880 \\
	11	   0.773228908695097 \\
	12	   0.019198423060312 \\
};
\addlegendentry{$\epsilon_{\text{state}}$}

\end{axis}
\end{tikzpicture}%
}
	\caption{Inline filter with finite transmission zeros error estimator results. (a) Convergence of the greedy algorithm. (b) Effectivity (\texttt{eff}).}
	\label{fig:Sec-NumericalResults-Subsec-MacchiarellaFilter-Error}
\end{figure}
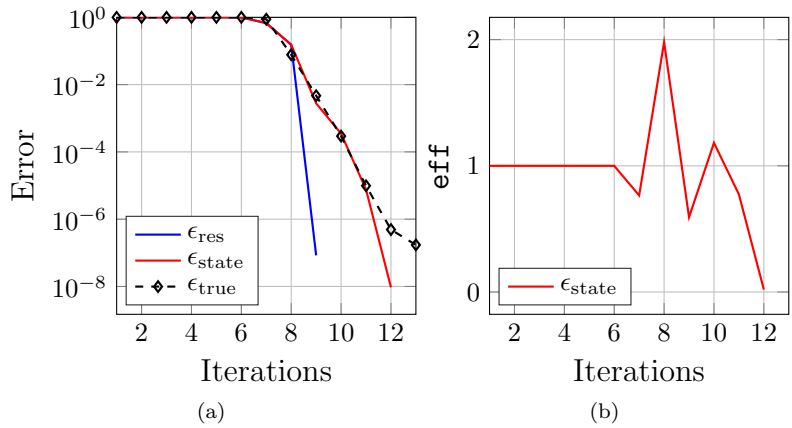

Next, we carry out the numerical tests to compare the different MOR methodologies. Table \ref{tab:Sec-NumericalResults-Subsec-MacchiarellaFilter-GreedySamples} details the frequency samples adaptively chosen by each greedy algorithm as well as the error estimator and true error at each iteration in the MOR process. Contrary to what Algorithm \ref{alg:Sec-InfSupConstantFree-SubSec-EnhancedReducedBasisSpace-FastROMErrorEstimator} proposes, it should be noted that the true error $\epsilon_{\text{true}}$-based greedy algorithm does not choose an end point sample after the eigenbasis made of 6 in-band resonant modes is built up. Fig. \ref{fig:Sec-NumericalResults-Subsec-MacchiarellaFilter-Error} depicts the convergence behavior for the different methodologies and details the effectivity metric for the proposed \emph{a posteriori} state error estimator. Reasonable performance is observed in the proposed methodology in comparison to the true error. Once again, the residual norm-based greedy algorithm prematurely aborts the iterative MOR procedure.

In addition, Table \ref{tab:Sec-NumericalResults-Subsec-MacchiarellaFilter-AlgorithmComparison} compares the performance of the proposed approach to Algorithm~\ref{alg:Sec-InfSupConstantFree-ROMErrorEstimator} where $\xi$ in Step \ref{step:Sec-InfSupConstantFree-Alg-ROMErrorEstimator-StoppingCriterion} is replaced by the indicator in \eqref{eq:Sec-NumericalResults-StateErrorEstimator}. Once again, even though the in-band eigenmodes are not \emph{a priori} included in the projection basis in Algorithm~\ref{alg:Sec-InfSupConstantFree-ROMErrorEstimator}, these show up in the greedy algorithm after the first few iterations. This shows Algorithm~\ref{alg:Sec-InfSupConstantFree-ROMErrorEstimator} is properly working out since the in-band eigenmodes are expected to show up in a good approximation basis, as has been discussed from a theoretical point of view. As far as computational burden is concerned, Algorithm~\ref{alg:Sec-InfSupConstantFree-ROMErrorEstimator} is more time consuming than Algorithm \ref{alg:Sec-InfSupConstantFree-SubSec-EnhancedReducedBasisSpace-FastROMErrorEstimator}. As a matter of fact, the size of both ROMs when the same stopping criterion is used is $12$ for Algorithm \ref{alg:Sec-InfSupConstantFree-SubSec-EnhancedReducedBasisSpace-FastROMErrorEstimator} and $14$ for Algorithm \ref{alg:Sec-InfSupConstantFree-ROMErrorEstimator}, whereas the size of the residual ROMs defined in  \eqref{eq:Sec-InfSupConstantFree-GalerkinProjectionError} is $13$ and $27$ for Algorithms  \ref{alg:Sec-InfSupConstantFree-SubSec-EnhancedReducedBasisSpace-FastROMErrorEstimator} and \ref{alg:Sec-InfSupConstantFree-ROMErrorEstimator}, respectively. This shows the additional effort that has to be carried out in Algorithm~\ref{alg:Sec-InfSupConstantFree-ROMErrorEstimator} with respect to Algorithm~\ref{alg:Sec-InfSupConstantFree-SubSec-EnhancedReducedBasisSpace-FastROMErrorEstimator}. As a result, the proposed approach gives rise to a fast \emph{a posteriori} state error estimation.

\begin{table}[tbp]
	\centering
	\begin{threeparttable}
		\renewcommand{\arraystretch}{1.3}
		\caption{Greedy Algorithm Frequency Samples for Different MOR Methodologies in the Inline Filter with Finite Transmission Zeros.}
		\label{tab:Sec-NumericalResults-Subsec-MacchiarellaFilter-GreedySamples}
		\centering
		\begin{tabular}{c|c|c|c|c|c}
			\hline
			GHz & $\epsilon_{\text{true}}$ & GHz & $\epsilon_{\text{state}}$ & GHz & $\epsilon_{\text{res}}$ \\ 
			\hline
			\hline
			1.6981     & 1. 				& 1.6981 	& 1. 				& 1.6981	& 1. \\
			1.7043     & 1. 				& 1.7043 	& 1. 				& 1.7043	& 1. \\
			1.7799     & 1. 				& 1.7799 	& 1. 				& 1.7799	& 1. \\
			1.8328     & 1. 				& 1.8328 	& 1. 				& 1.8328	& 1. \\
			2.7661     & 1. 				& 2.7661 	& 1. 				& 2.7661	& 1. \\
			2.7856     & 1. 				& 2.7856 	& 1. 				& 2.7856	& 1. \\
			2.3000     & $8.8\cdot10^{-1}$ 	& 1.4000   	& $6.7\cdot10^{-1}$ & 1.4000   	& $6.7\cdot10^{-1}$ \\
			2.9000     & $7.7\cdot10^{-2}$ 	& 2.9000	& $1.5\cdot10^{-1}$ & 2.9000   	& $1.5\cdot10^{-1}$ \\
			1.4000     & $4.6\cdot10^{-3}$ 	& 2.5156   	& $2.8\cdot10^{-3}$ & 2.7661    & $8.4\cdot10^{-8}$ \\
			2.6200     & $2.9\cdot10^{-4}$ 	& 2.0520	& $3.5\cdot10^{-4}$ & --    	& -- \\
			1.9800     & $9.8\cdot10^{-6}$	& 2.8187	& $7.6\cdot10^{-6}$	& --    	& -- \\
			2.8500     & $4.9\cdot10^{-7}$	& 1.7007	& $9.4\cdot10^{-9}$	& --    	& -- \\
			1.7600     & $1.7\cdot10^{-7}$	& --		& -- 				& --    	& -- \\
			\hline
		\end{tabular}
	\end{threeparttable}
\end{table}
\begin{table}[tbp]
	\centering
	\begin{threeparttable}
		\renewcommand{\arraystretch}{1.3}
		\caption{State Error Algorithm Comparison in the Inline Filter with Finite Transmission Zeros.}
		\label{tab:Sec-NumericalResults-Subsec-MacchiarellaFilter-AlgorithmComparison}
		\centering
		\begin{tabular}{c|c|c|c}
			\hline
			\multicolumn {2} {c|}{Algorithm \ref{alg:Sec-InfSupConstantFree-SubSec-EnhancedReducedBasisSpace-FastROMErrorEstimator}} & \multicolumn {2} {|c}{Algorithm \ref{alg:Sec-InfSupConstantFree-ROMErrorEstimator}} \\ 
			\hline
			\hline
			GHz & $\epsilon_{\text{state}}$ & GHz & $\epsilon_{\text{state}}$ \\ 
			\hline
			\hline
			1.6981 	& 1. 				& 1.4000	& $3.6\cdot10^{-3}$ \\
			1.7043 	& 1. 				& 1.9430	& $1.6\cdot10^{-2}$ \\
			1.7799 	& 1. 				& 1.7260	& $6.1$ \\
			1.8328 	& 1. 				& 2.7660	& $1.2\cdot10^{-2}$ \\
			2.7661 	& 1. 				& 1.6980	& $6.3\cdot10^{-4}$ \\
			2.7856 	& 1. 				& 2.7860	& $2.5\cdot10^{-6}$ \\
			1.4000  & $6.7\cdot10^{-1}$ & 1.8330  	& $3.6\cdot10^{-8}$ \\
			2.9000	& $1.5\cdot10^{-1}$ & --	  	& -- \\
			2.5156  & $2.8\cdot10^{-3}$ & --   		& -- \\
			2.0520	& $3.5\cdot10^{-4}$ & --    	& -- \\
			2.8187	& $7.6\cdot10^{-6}$	& --    	& -- \\
			1.7007	& $9.4\cdot10^{-9}$	& --    	& -- \\			
			\hline
			\hline
			ROM size 			& 12 	& ROM size			& 14 \\
			Residual ROM size 	& 13 	& Residual ROM size	& 27 \\
			\hline
		\end{tabular}
	\end{threeparttable}
\end{table}
\subsection{Inline Dielectric Resonator Filter}
\label{Sec-NumericalResults-Subsec-SnyderFilter}
\begin{figure}[tbp]
	\centering
	\includegraphics[width=\linewidth]{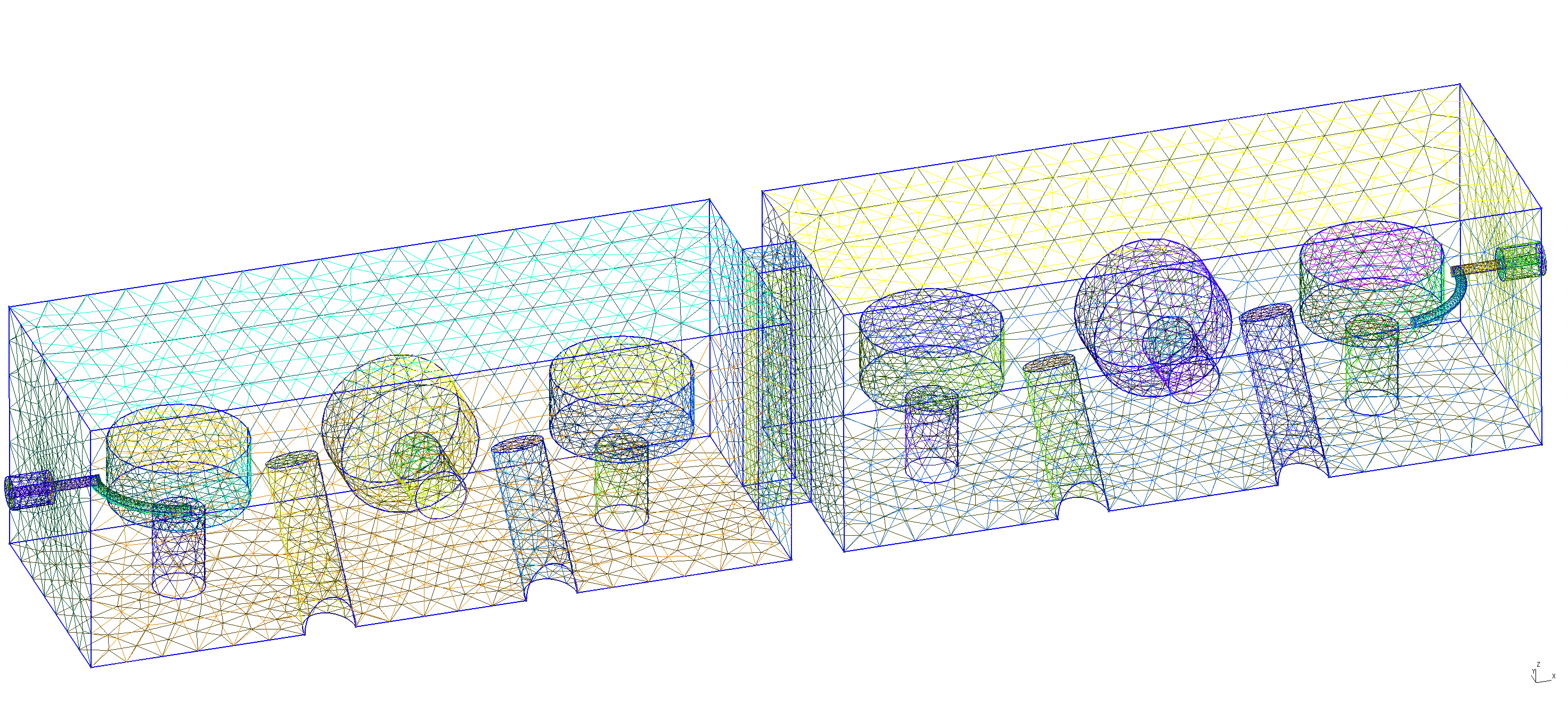}
	\caption{Inline dielectric resonator filter designed in \cite{bastioli2012inline}.}
	\label{fig:Sec-NumericalResults-Subsec-SnyderFilter-FilterGeometry}
\end{figure}

\begin{figure}[tbp]
	\centering
	\input{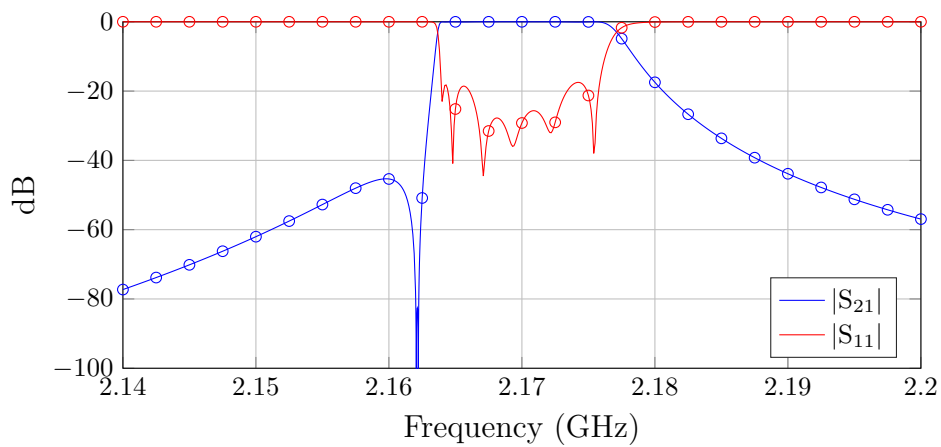}
	\caption{Inline dielectric resonator filter scattering parameter response comparison with the proposed approach. (--) RBM. ($\circ$) FEM.}
	\label{fig:Sec-NumericalResults-Subsec-SnyderFilter-FilterResponse}
\end{figure}

A sixth order inline dielectric resonator filter with two transmission zeros is depicted in Fig. \ref{fig:Sec-NumericalResults-Subsec-SnyderFilter-FilterGeometry}. Cross-coupling between nonadjacent dielectric resonators, appropriately arranging their orientations, is obtained by exploiting multiple evanescent modes in the inline structure. This filter is proposed in \cite{bastioli2012inline}. The $[2.14, 2.2]$ GHz band is taken into account in the analysis. A FEM discretization shown in Fig. \ref{fig:Sec-NumericalResults-Subsec-SnyderFilter-FilterGeometry} with 230,058 degrees of freedom is used. Fig. \ref{fig:Sec-NumericalResults-Subsec-SnyderFilter-FilterResponse} details the filter response in the band of analysis under both time-consuming FEM simulation and fast RBM analysis. A ROM of dimension 10 is used to get the fast frequency sweep results in Fig. \ref{fig:Sec-NumericalResults-Subsec-SnyderFilter-FilterResponse}. Good agreement is obtained between both analyses. It should be pointed out that the FEM solution of the FOM evaluated at a given frequency takes $5.610$ seconds. In the online stage, the ROM resulting from RBM was evaluated for $1201$ different frequency samples requiring $0.085$ seconds in total; this works out to $70$ microseconds to solve a single ROM and amounts to a speedup of nearly $80000$.

\begin{figure}[tbp]
	\centering
	\subfloat[]{\label{fig:Sec-NumericalResults-Subsec-SnyderFilter-ErrorConvergence}
%
%
\definecolor{mycolor1}{rgb}{0.00000,0.44700,0.74100}%
\begin{tikzpicture}

\begin{axis}[%
width=0.3\fwidth,
height=0.3\fheight,
at={(0\fwidth,0\fheight)},
scale only axis,
xmin=1,
xmax=11,
xtick = {2,4,6,8,10},
xlabel style={font=\color{white!15!black}, font = \large},
xlabel={Iterations},
ymode=log,
ymax=1e0,
yminorticks=true,
ytick = {1e0, 1e-2, 1e-4, 1e-6, 1e-8 },
ylabel style={font=\color{white!15!black}, font = \large},
ylabel={Error},
axis background/.style={fill=white},
xmajorgrids,
ymajorgrids,
yminorgrids,
legend pos=south west,
legend style={legend cell align=left, align=left, draw=white!15!black}
]
\addplot [color=blue, thick]
  table[row sep=crcr]{%
1	1.0\\
2	1.0\\
3	1.0\\
4	1.0\\
5	1.0\\
6	1.0\\
7   0.341493777395724\\
8   0.000000054220845\\
};
\addlegendentry{$\epsilon_{\text{res}}$}

\addplot [color=red, thick]
table[row sep=crcr]{%
1	1.0\\
2	1.0\\
3	1.0\\
4	1.0\\
5	1.0\\
6	1.0\\
7   0.341493777395724\\
8   0.007828741916809\\
9   0.000005686000352\\
10  0.000000166281689\\
};
\addlegendentry{$\epsilon_{\text{state}}$}

\addplot [color=black, dashed, mark=diamond, mark options={solid, black}, thick]
  table[row sep=crcr]{%
1	1.0\\
2	1.0\\
3	1.0\\
4	1.0\\
5	1.0\\
6	1.0\\
7   0.341493923811244\\
8   0.007828744471497\\
9   0.000027630282300\\
10  0.000000364705772\\
11  0.000000161202543\\
};
\addlegendentry{$\epsilon_{\text{true}}$}

\end{axis}
\end{tikzpicture}%
} 
	\subfloat[]{\label{fig:Sec-NumericalResults-Subsec-SnyderFilter-ErrorEffectivity}
%
%
\definecolor{mycolor1}{rgb}{0.00000,0.44700,0.74100}%
\begin{tikzpicture}

\begin{axis}[%
width=0.3\fwidth,
height=0.3\fheight,
at={(0\fwidth,0\fheight)},
scale only axis,
xmin=1,
xmax=11,
xtick = {2,4,6,8,10,12},
xlabel style={font=\color{white!15!black}, font = \large},
xlabel={Iterations},
yminorticks=true,
ylabel style={font=\color{white!15!black}, font = \large},
ylabel={\texttt{eff}},
axis background/.style={fill=white},
xmajorgrids,
ymajorgrids,
yminorgrids,
legend pos=south west,
legend style={legend cell align=left, align=left, draw=white!15!black}
]
\addplot [color=red, thick]
table[row sep=crcr]{%
1	   1.000000000000000\\
2	   1.000000000000000\\
3	   1.000000000000000\\
4	   1.000000000000000\\
5	   1.000000000000000\\
6	   1.000000000000000\\
7  	   0.999999571249997\\
8	   0.999999673678393\\
9      0.205788717245862\\
10	   0.455933800769630\\
};
\addlegendentry{$\epsilon_{\text{state}}$}

\end{axis}
\end{tikzpicture}%
}
	\caption{Inline dielectric resonator filter error estimator results. (a) Convergence of the greedy algorithm. (b) Effectivity (\texttt{eff}).}
	\label{fig:Sec-NumericalResults-Subsec-SnyderFilter-Error}
\end{figure}
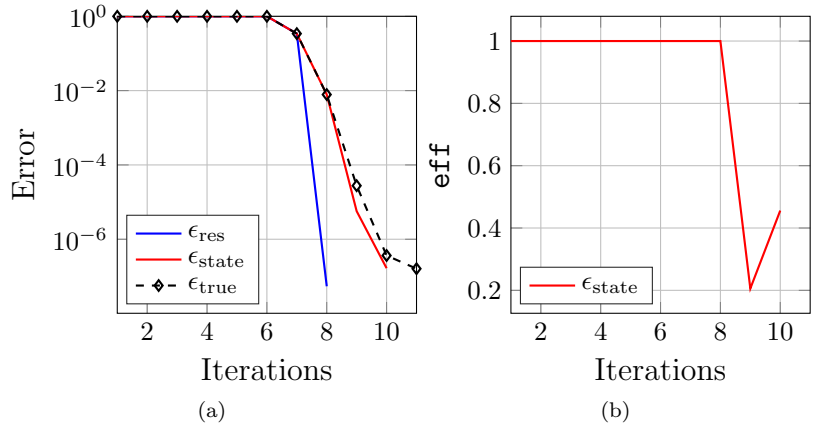

Next, we compare the different MOR techniques. Table~\ref{tab:Sec-NumericalResults-Subsec-SnyderFilter-GreedySamples} shows the frequency samples adaptively chosen by each greedy algorithm as well as the error estimator at each iteration in the MOR process. As expected in a sixth order filter, 6 in-band eigenmodes are found in the frequency band of analysis. Fig. \ref{fig:Sec-NumericalResults-Subsec-SnyderFilter-Error} details the convergence behavior for the different methodologies and shows the effectivity metric for the proposed approach. A good behavior is observed in the proposed \emph{a posteriori} state error estimator.
\begin{table}[tbp]
	\centering
	\begin{threeparttable}
		\renewcommand{\arraystretch}{1.3}
		\caption{Greedy Algorithm Frequency Samples for Different MOR Methodologies in the Inline Dielectric Resonator Filter.}
		\label{tab:Sec-NumericalResults-Subsec-SnyderFilter-GreedySamples}
		\centering
		\begin{tabular}{c|c|c|c|c|c}
			\hline
			GHz & $\epsilon_{\text{true}}$ & GHz & $\epsilon_{\text{state}}$ & GHz & $\epsilon_{\text{res}}$ \\ 
			\hline
			\hline
			2.1635     & 1. 				& 2.1635 	& 1. 				& 2.1635	& 1. \\
			2.1640     & 1. 				& 2.1640 	& 1. 				& 2.1640	& 1. \\
			2.1663     & 1. 				& 2.1663 	& 1. 				& 2.1663	& 1. \\
			2.1709     & 1. 				& 2.1709 	& 1. 				& 2.1709	& 1. \\
			2.1768     & 1. 				& 2.1768 	& 1. 				& 2.1768	& 1. \\
			2.1788     & 1. 				& 2.1788 	& 1. 				& 2.1788	& 1. \\
			2.1400     & $3.4\cdot10^{-1}$ 	& 2.1400   	& $3.4\cdot10^{-1}$ & 2.1400   	& $3.4\cdot10^{-1}$ \\			
			2.2000     & $7.8\cdot10^{-3}$ 	& 2.2000	& $7.8\cdot10^{-3}$ & 2.1709   	& $5.4\cdot10^{-8}$ \\			
			2.1530     & $2.8\cdot10^{-5}$	& 2.1697	& $5.7\cdot10^{-6}$ & --      	& -- \\			
			2.1915     & $3.6\cdot10^{-7}$	& 2.1874	& $1.7\cdot10^{-7}$ & --      	& -- \\			
			2.1785     & $1.6\cdot10^{-7}$	& --		& -- 				& --      	& -- \\			
			\hline
		\end{tabular}
	\end{threeparttable}
\end{table}
\subsection{Combline Diplexer}
\label{Sec-NumericalResults-Subsec-ComblineDiplexer}
\begin{figure}[tbp]
	\centering
	\includegraphics[width=\linewidth]{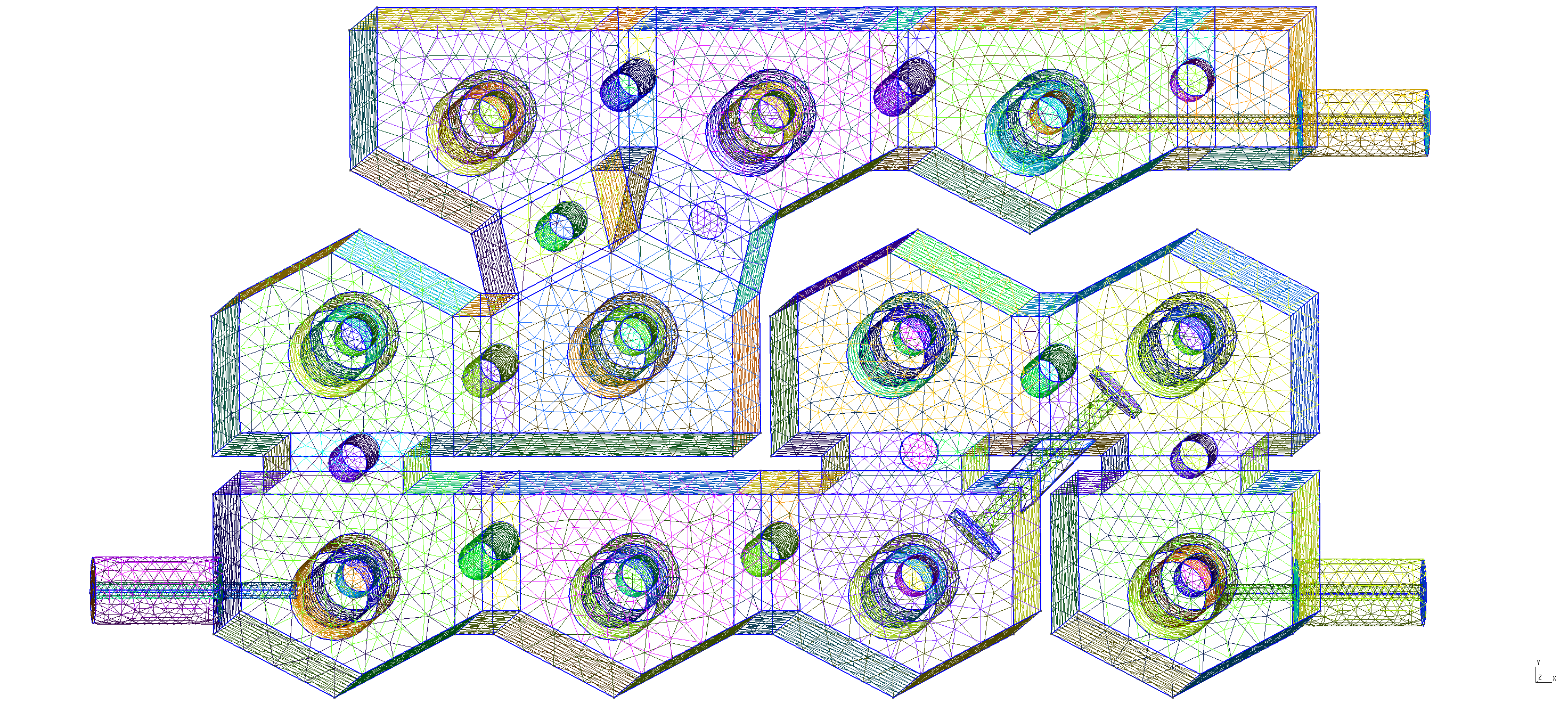}
	\caption{Combline diplexer designed in \cite{zhao2014iterative}.}
	\label{fig:Sec-NumericalResults-Subsec-ComblineDiplexer-FilterGeometry}
\end{figure}

The last real-life application is an $11^\text{th}$ order combline diplexer with star-junction designed in \cite{zhao2014iterative}. The geometry of this diplexer is shown in Fig. \ref{fig:Sec-NumericalResults-Subsec-ComblineDiplexer-FilterGeometry}. The frequency band of analysis is $\mathcal{B}:=[2.2, 3.0]$ GHz. An FEM system with 270,446 degrees of freedom arises whereas the application of the proposed methodology in Algorithm \ref{alg:Sec-InfSupConstantFree-SubSec-EnhancedReducedBasisSpace-FastROMErrorEstimator} gives rise to a ROM of dimension 20 by means of RBM. The scattering parameter response for this diplexer is detailed in Fig. \ref{fig:Sec-NumericalResults-Subsec-ComblineDiplexer-FilterResponse}. Good agreement is found between FEM and RBM results. It should be pointed out that further tuning is needed to obtain the target equiripple response. The FOM solution for this example requires $4.630$ seconds. In the online stage, the ROM was evaluated at $1601$ different frequencies taking $0.446$ seconds. Thus, the time to solve a single ROM is $278$ microseconds, at a speedup of around $16650$.

As expected in an $11^\text{th}$ order diplexer, 11 in-band eigenmodes are found in the frequency band of interest. Table \ref{tab:Sec-NumericalResults-Subsec-ComblineDiplexer-GreedySamples} details the different frequency samples for each methodology. A comparison for the different MOR strategies is shown in Fig. \ref{fig:Sec-NumericalResults-Subsec-ComblineDiplexer-Error}, where the convergence of the estimated and true errors at each iteration as well as the effectivity metric is detailed. Once again, the residual norm-based greedy algorithm prematurely stops, misled by oversampling nearby the eigenresonances. A reasonable performance is observed in the proposed \emph{a posteriori} state error estimator.
\begin{table}[tbp]
	\centering
	\begin{threeparttable}
		\renewcommand{\arraystretch}{1.3}
		\caption{Greedy Algorithm Frequency Samples for Different MOR Methodologies in the Combline Diplexer.}
		\label{tab:Sec-NumericalResults-Subsec-ComblineDiplexer-GreedySamples}
		\centering
		\begin{tabular}{c|c|c|c|c|c}
			\hline
			GHz & $\epsilon_{\text{true}}$ & GHz & $\epsilon_{\text{state}}$ & GHz & $\epsilon_{\text{res}}$ \\ 
			\hline
			\hline
   2.4623 	&	1.				&   2.4623 	&	1.				&  	2.4623 	&	1.	\\
   2.4914	&	1.				&   2.4914	&	1.				&   2.4914	&	1.	\\
   2.5233	&	1.				&   2.5233	&	1.				&   2.5233	&	1.	\\
   2.5592	&	1.				&   2.5592	&	1.				&   2.5592	&	1.	\\
   2.5803	&	1.				&   2.5803	&	1.				&   2.5803	&	1.	\\
   2.5981	&	1.				&   2.5981	&	1.				&   2.5981	&	1.	\\
   2.6163	&	1.				&   2.6163	&	1.				&   2.6163	&	1.	\\
   2.6419	&	1.				&   2.6419	&	1.				&   2.6419	&	1.	\\
   2.6766	&	1.				&   2.6766	&	1.				&   2.6766	&	1.	\\
   2.7102	&	1.				&   2.7102	&	1.				&   2.7102	&	1.	\\
   2.7274	&	1.				&   2.7274	&	1.				&   2.7274	&	1.	\\
   3.0000	& $8.6\cdot10^{-1}$ &   3.0000	& $8.6\cdot10^{-1}$ &	3.0000  &	$8.6\cdot10^{-1}$ \\
   2.2000	& $4.5\cdot10^{-1}$ &   2.2000	& $4.5\cdot10^{-1}$ &	2.2000  &	$4.5\cdot10^{-1}$ \\
   2.8600	& $6.0\cdot10^{-2}$ &   2.5094	& $1.1\cdot10^{-2}$ &	2.6766  &	$3.1\cdot10^{-9}$ \\
   2.3500	& $7.2\cdot10^{-3}$ &	2.9067	& $1.1\cdot10^{-2}$ &	-- 		&	-- \\
   2.9500	& $1.1\cdot10^{-3}$ &	2.6060	& $9.2\cdot10^{-5}$ &	-- 		&	-- \\
   2.2600	& $1.6\cdot10^{-4}$ &	2.2835	& $2.4\cdot10^{-4}$ &	-- 		&	-- \\
   2.7800	& $1.5\cdot10^{-5}$&	2.9712	& $3.7\cdot10^{-5}$ &	-- 		&	-- \\
   2.4200	& $1.1\cdot10^{-6}$&	2.2137	& $7.9\cdot10^{-7}$ &	-- 		&	-- \\
   2.4100	& $1.9\cdot10^{-7}$&	2.6766	& $1.1\cdot10^{-8}$ &	-- 		&	-- \\   
			\hline
		\end{tabular}
	\end{threeparttable}
\end{table}
\begin{figure}[tbp]
	\centering
	\input{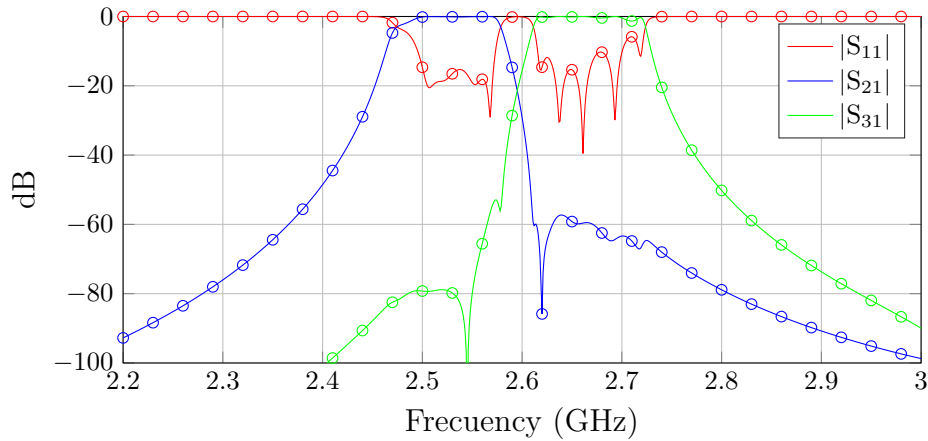}
	\caption{Combline diplexer scattering parameter response comparison with the proposed approach. (--)~RBM. ($\circ$)~FEM.}
	\label{fig:Sec-NumericalResults-Subsec-ComblineDiplexer-FilterResponse}
\end{figure}	
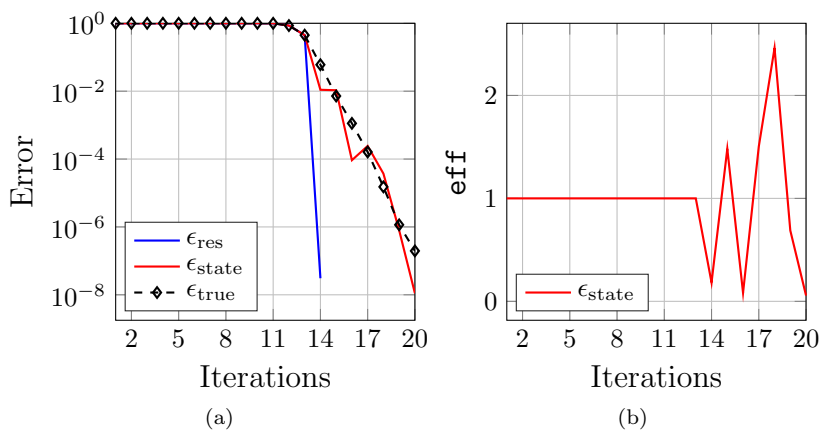
\begin{figure}[tbp]
	\centering
	\subfloat[]{\label{fig:Sec-NumericalResults-Subsec-ComblineDiplexer-ErrorConvergence}
%
%
\definecolor{mycolor1}{rgb}{0.00000,0.44700,0.74100}%
\begin{tikzpicture}

\begin{axis}[%
width=0.3\fwidth,
height=0.3\fheight,
at={(0\fwidth,0\fheight)},
scale only axis,
xmin=1,
xmax=20,
xtick = {2,5,8,11,14,17,20},
xlabel style={font=\color{white!15!black}, font = \large},
xlabel={Iterations},
ymode=log,
ymax=1e0,
yminorticks=true,
ytick = {1e0, 1e-2, 1e-4, 1e-6, 1e-8},
ylabel style={font=\color{white!15!black}, font = \large},
ylabel={Error},
axis background/.style={fill=white},
xmajorgrids,
ymajorgrids,
yminorgrids,
legend pos=south west,
legend style={legend cell align=left, align=left, draw=white!15!black}
]
\addplot [color=blue, thick]
  table[row sep=crcr]{%
1	1.0 \\
2	1.0 \\
3	1.0 \\
4	1.0 \\
5	1.0 \\
6	1.0 \\
7	1.0 \\
8	1.0 \\
9	1.0 \\
10	1.0 \\
11	1.0 \\
12   0.857011668531999 \\
13   0.447458377952631 \\
14   0.000000030786572 \\
};
\addlegendentry{$\epsilon_{\text{res}}$}

\addplot [color=red, thick]
table[row sep=crcr]{%
1	1.0 \\
2	1.0 \\
3	1.0 \\
4	1.0 \\
5	1.0 \\
6	1.0 \\
7	1.0 \\
8	1.0 \\
9	1.0 \\
10	1.0 \\
11	1.0 \\
12   0.857011668531999\\
13   0.447458377952631\\
14   0.010979070998951\\
15   0.010666770832825\\
16   0.000092606047319\\
17   0.000243163730848\\
18   0.000037275461097\\
19   0.000000793982368\\
20   0.000000011098919\\
};
\addlegendentry{$\epsilon_{\text{state}}$}

\addplot [color=black, dashed, mark=diamond, mark options={solid, black}, thick]
  table[row sep=crcr]{%
1	1.0\\
2	1.0\\
3	1.0\\
4	1.0\\
5	1.0\\
6	1.0\\
7	1.0\\
8	1.0\\
9	1.0\\
10	1.0\\
11	1.0\\
12   0.857011668531999 \\
13   0.447458713179216 \\
14   0.060199435213297 \\
15   0.007196918090405 \\
16   0.001125739756782 \\
17   0.000161669601348 \\
18   0.000015150696354 \\
19   0.000001153728304 \\
20   0.000000196497992 \\
};
\addlegendentry{$\epsilon_{\text{true}}$}

\end{axis}
\end{tikzpicture}%
} 
	\subfloat[]{\label{fig:Sec-NumericalResults-Subsec-ComblineDiplexer-ErrorEffectivity}
%
%
\definecolor{mycolor1}{rgb}{0.00000,0.44700,0.74100}%
\begin{tikzpicture}

\begin{axis}[%
width=0.3\fwidth,
height=0.3\fheight,
at={(0\fwidth,0\fheight)},
scale only axis,
xmin=1,
xmax=20,
xtick = {2,5,8,11,14,17,20},
xlabel style={font=\color{white!15!black}, font = \large},
xlabel={Iterations},
yminorticks=true,
ylabel style={font=\color{white!15!black}, font = \large},
ylabel={\texttt{eff}},
axis background/.style={fill=white},
xmajorgrids,
ymajorgrids,
yminorgrids,
legend pos=south west,
legend style={legend cell align=left, align=left, draw=white!15!black}
]
\addplot [color=red, thick]
table[row sep=crcr]{%
1	1.000000000000000 \\
2	1.000000000000000 \\
3	1.000000000000000 \\
4	1.000000000000000 \\
5	1.000000000000000 \\
6	1.000000000000000 \\
7	1.000000000000000 \\
8	1.000000000000000 \\
9	1.000000000000000 \\
10	1.000000000000000 \\
11	1.000000000000000 \\
12	1.000000000000000 \\
13  0.999999250821194 \\
14  0.182378305710840 \\
15  1.482130364530010 \\
16  0.082262393915483 \\
17  1.504078248605605 \\
18  2.460313389367199 \\
19  0.688188340916860 \\
20  0.05648362554259588 \\
};
\addlegendentry{$\epsilon_{\text{state}}$}

\end{axis}
\end{tikzpicture}%
}
	\caption{Combline diplexer error estimator results. (a) Convergence of the greedy algorithm. (b) Effectivity (\texttt{eff}).}
	\label{fig:Sec-NumericalResults-Subsec-ComblineDiplexer-Error}
\end{figure}	
\section{Conclusions}
\label{Sec-Conclusions}
A compact and reliable MOR method for fast frequency sweeps in microwave circuits by means of the reduced-basis method has been detailed. A compact basis including the in-band resonant modes hit in the frequency band of interest for both the reduced basis approximation of the electric field and its corresponding state error has been proposed. This allows to efficiently solve for both the reduced original and residual problems, thus minimizing the additional computational effort. The benefits of using proper state error estimators avoiding time-consuming \emph{inf-sup} constant evaluations has also been highlighted. As a result, a fast \emph{a posteriori} state error estimator for the ROM has been obtained. Real-life microwave devices, including a quad-mode dielectric resonator filter and a combline diplexer, have shown the capabilities and reliability of the proposed methodology.
%

\section*{Acknowledgments}%
\addcontentsline{toc}{section}{Acknowledgments}

Sridhar Chellappa is supported by the International Max Planck Research School for Advanced Methods in Process and Systems Engineering (IMPRS-ProEng).


\addcontentsline{toc}{section}{References}
\bibliographystyle{plainurl}
\bibliography{references}

\end{document}